\documentclass[12pt]{article}
\usepackage{ejc}
\usepackage[reqno]{amsmath}
\usepackage{amssymb, theorem, enumerate, bm}
\usepackage{url}
\usepackage{hyperref}
\hypersetup{colorlinks=false}

\newcommand{\rank}{     \mathrm{rank}\, }

\newcommand{\nul}{      \mathrm{nullity}\,  }
\newcommand{\ch}{      \mathrm{char}\,  }

\newcommand\cF{{\mathcal F}}
\newcommand\cG{{\mathcal G}}

\newcommand\cR{{\mathcal R}}
\newcommand\cS{{\mathcal S}}

\makeatletter
\g@addto@macro\bfseries{\boldmath}
\makeatother

\expandafter\ifx\csname pdfoptionalwaysusepdfpagebox\endcsname\relax\else
\fi

\theoremstyle{plain}
{\theorembodyfont{\slshape}
\newtheorem{theorem}{Theorem}[section]
\newtheorem{lemma}[theorem]{Lemma}
\newtheorem{corollary}[theorem]{Corollary}}
\theorembodyfont{\rmfamily}

\newtheorem{definition}[theorem]{Definition}

\newtheorem{conjecture}[theorem]{Conjecture}
\newtheorem{question}[theorem]{Question}

\newcommand\lref[1]{Lemma~\ref{lem:#1}}
\newcommand\tref[1]{Theorem~\ref{thm:#1}}
\newcommand\cref[1]{Corollary~\ref{cor:#1}}

\newcommand\sref[1]{Section~\ref{sec:#1}}
\newcommand\conjref[1]{Conjecture~\ref{conj:#1}}
\newcommand\dref[1]{Definition~\ref{def:#1}}

\newcommand\qref[1]{Question~\ref{ques:#1}}

\def\sqr#1#2{{\vbox{\hrule height.#2pt
    \hbox{\vrule width.#2pt height#1pt \kern#1pt
        \vrule width.#2pt}\hrule height.#2pt}}}
\def\eqed{\sqr53}
\def\qed{%
    \ifmmode\eqno\eqed
    \else\nobreak\ \hfill\eqed\medbreak\fi}

\newcommand{\upperRomannumeral}[1]{\uppercase\expandafter{\romannumeral#1}}

\title{Inclusion Matrices and the MDS Conjecture}
\author{Ameera Chowdhury \thanks{Department of Mathematics, Rutgers University, Piscataway, NJ, 08854-8019, USA. E-mail: {\tt ameerah@alumni.caltech.edu.} Research supported by NSF grant DMS-1203982.}}

\begin{document}
\maketitle

\begin{abstract}
Let $\mathbb{F}_{q}$ be a finite field of order $q$ with characteristic $p$. An arc in $\mathbb{F}_{q}^{k}$ is an ordered family of at least $k$ vectors in which every subfamily of size $k$ is a basis of $\mathbb{F}_{q}^{k}$. The MDS conjecture, which was posed by Segre in 1955, states that if $k \leq q$, then an arc in $\mathbb{F}_{q}^{k}$ has size at most $q+1$, unless $q$ is even and $k=3$ or $k=q-1$, in which case it has size at most $q+2$. 

We propose a conjecture which would imply that the MDS conjecture is true for almost all values of $k$ when $q$ is odd. We prove our conjecture in two cases and thus give simpler proofs of the MDS conjecture when $k \leq p$, and if $q$ is not prime, for $k \leq 2p-2$. To accomplish this, given an arc $G \subset \mathbb{F}_{q}^{k}$ and a nonnegative integer $n$, we construct a matrix $M_{G}^{\uparrow n}$, which is related to an inclusion matrix, a well-studied object in combinatorics. Our main results relate algebraic properties of the matrix $M_{G}^{\uparrow n}$ to properties of the arc $G$ and may provide new tools in the computational classification of large arcs.
\end{abstract}

\section{Introduction}
\label{sec:intro}

Let $\mathbb{F}_{q}$ be a finite field of order $q$ with characteristic $p$. An \textsl{arc} in $\mathbb{F}_{q}^{k}$ is an ordered family of at least $k$ vectors in which every subfamily of size $k$ is a basis of $\mathbb{F}_{q}^{k}$. Most authors define an arc, equivalently, as an unordered set of points in the corresponding projective space. For the techniques developed in this article, however, we find it more convenient to define arcs as ordered families of vectors. On the other hand, we will denote arcs with set notation rather than tuple notation as this is more natural.

Given an arc $G \subset \mathbb{F}_{q}^{k}$ and a basis $B$ of $\mathbb{F}_{q}^{k}$, let $M(G,B)$ be the matrix whose columns are the vectors in $G$ written with respect to the basis $B$ in the order given by $G$. If $G, G' \subset \mathbb{F}_{q}^{k}$ are two arcs, then we say that $G$ is \textsl{linearly equivalent} to $G'$ if the matrix $M(G,B)$ can be transformed into the matrix $M(G',B)$ using only elementary row operations, column permutations, and multiplication of columns by nonzero scalars.

A natural question is to determine how large an arc in $\mathbb{F}_{q}^{k}$ can be.

\begin{question}
\label{ques:motivation}
What is the maximum size $g(k,q)$ of an arc in $\mathbb{F}_{q}^{k}$?
\end{question}

\noindent \qref{motivation} interests the coding theory, algebraic geometry, and finite geometry communities, and its importance is highlighted by a \$1000 prize offered for its solution by the Information Theory and Applications (ITA) center at UCSD \cite{ITAPrize}.

If $(e_{1}, \ldots, e_{k})$ is a basis for $\mathbb{F}_{q}^{k}$, then a natural arc in $\mathbb{F}_{q}^{k}$ of size $k+1$ is given by
\begin{equation}
\label{bushexample}
\{e_{1}, \ldots, e_{k}, e_{1}+\cdots+e_{k}\},
\end{equation}
which proves that $g(k,q) \geq k+1$. A straightforward argument shows that $g(k,q)=k+1$ when $k \geq q$, and moreover if $S \subset \mathbb{F}_{q}^{k}$ is an arc of size $k+1$, then $S$ is linearly equivalent to \eqref{bushexample}. This result was first proved by Bush \cite{Bush} in 1952.

\qref{motivation} becomes difficult to answer, however, when $k < q$. In this case, we can construct arcs that are larger than the arc in \eqref{bushexample}. For example, the \textsl{normal rational curve} $\cR_{k} \subset \mathbb{F}_{q}^{k}$, which is defined by
\begin{equation}
\label{vandermondeexample}
\cR_{k} = \{(1,t,t^2, \ldots, t^{k-1}) \, | \, t \in \mathbb{F}_{q}\} \cup \{(0, \ldots, 0, 1)\},
\end{equation}
is an arc of size $q+1$. The normal rational curve $\cR_{k}$ shows that $g(k,q) \geq q+1$, and in 1955, Segre \cite{MDSConj} conjectured that this lower bound is tight in most cases when $k \leq q$.

\begin{conjecture}[Segre, \cite{MDSConj}]
\label{conj:main}
If $k \leq q$, then the maximum size $g(k,q)$ of an arc in $\mathbb{F}_{q}^{k}$ is
\[
g(k,q) = \begin{cases}
         q+1 & \mbox{if $q$ is odd or $k \notin \{3,q-1\}$} \\
         q+2 & \mbox{if $q$ is even and $k \in \{3,q-1\}$}.
         \end{cases}
\]         
\end{conjecture}
\noindent \conjref{main} is called the MDS conjecture or the main conjecture for maximum distance separable codes, and was first posed by Segre as a question.

By the well-known principle of duality, if $S \subset \mathbb{F}_{q}^{k}$ is an arc of size $s > k$, then up to linear equivalence, we can associate a unique dual arc $S^{\perp} \subset \mathbb{F}_{q}^{s-k}$ of size $s$. This has two immediate implications. First, it explains why in \conjref{main}, exceptions occur for both $k=3$ and $k=q-1$ when $q$ is even. Second, it shows that if $g(k,q) = q+1$, then $g(q+2-k,q)=q+1$. As a result, if $q$ is odd and $g(k,q) = q+1$ when $k \leq (q+2)/2$, then $g(k,q) = q+1$ for all $k \leq q$. Duality thus allows us to prove \conjref{main} when $q$ is odd by restricting to the case $k \leq (q+2)/2$.

Ball \cite{MDSprime} proved that $g(k,q) = q+1$ when $k \leq p = \ch(\mathbb{F}_{q})$, and thus verified \conjref{main} when $q$ is prime. For a complete list of when \conjref{main} is known to hold for $q$ non-prime, see \cite{PackingSurvey} and \cite{GGG}. The best-known bounds up to first-order of magnitude ($c_{i}$ are constants), are that for $q$ an odd non-square, we have $g(k,q)=q+1$ when $k < \sqrt{pq}/4 + c_{1}p$, which was proved by Voloch \cite{Voloch}. For $q=p^{2h}$, where $p \geq 5$ is a prime, we have $g(k,q)=q+1$ when $k \leq \sqrt{q}/2 + c_{2}$, which was proved by Hirschfeld and Korchm\'{a}ros \cite{HK}. Ball and De Buele \cite{MDS2pminus2} proved that $g(k,q)=q+1$ when $k \leq 2\sqrt{q}-2$ and $q = p^{2}$.

If $k \leq q$ and $q$ is odd or $k \notin \{3,q-1\}$, it is natural to ask if the normal rational curve $\cR_{k}$ is the unique arc in $\mathbb{F}_{q}^{k}$ of size $q+1$ up to linear equivalence. By results of Kaneta and Maruta \cite{KanetaMaruta} and Seroussi and Roth \cite{SeroussiRoth}, a positive answer to this question would imply \conjref{main}. For many values of $k$ and $q$, the normal rational curve $\cR_{k}$ is the unique arc in $\mathbb{F}_{q}^{k}$ of size $q+1$ up to linear equivalence \cite{PackingSurvey}, but Glynn \cite{Glynn} showed that this is not always true. The Glynn arc $\cG \subset \mathbb{F}_{9}^{5}$ is an arc of size 10 and is defined by 
\begin{equation}
\label{glynn}
\cG = \{(1,t,t^{2} + \eta t^6, t^3, t^4) \, | \, t \in \mathbb{F}_{9}\} \cup \{(0,0,0,0,1)\},
\end{equation}
where $\eta \in \mathbb{F}_{9}$ satisfies $\eta^{4} = -1$. Remarkably, the Glynn arc $\cG$ is the only known arc in $\mathbb{F}_{q}^{k}$ of size $q+1$ that is not linearly equivalent to the normal rational curve $\cR_{k}$ when $k \leq q$ and $q$ is odd. 

\subsection{New Results}
\label{sec:ourresults}

We propose a conjecture, \conjref{PErrank}, which would imply that $g(k,q)=q+1$ when
\begin{equation}
\label{MDSbound}
k \leq \left( \frac{p-2}{2p-3} \right)q + \left(3 - \frac{p-1}{2p-3} \right),
\end{equation}
where $p = \ch(\mathbb{F}_{q})$. In \sref{intro}, we noted that to prove \conjref{main} when $q$ is odd, it suffices to restrict to the case $k \leq (q+2)/2$ by duality. As $p$ grows, the right hand side of \eqref{MDSbound} becomes very close to $(q+2)/2$. Consequently, if \conjref{PErrank} is true, then \conjref{main} is true for almost all values of $k$ when $q$ is odd.

To state \conjref{PErrank}, given an arc $G \subset \mathbb{F}_{q}^{k}$ and a nonnegative integer $n$, we define a matrix $M_{G}^{\uparrow n}$ whose algebraic properties are related to properties of $G$.

\begin{definition} 
\label{def:M}
Let $G \subset \mathbb{F}_{q}^{k}$ be an arc and let $0 \leq n \leq |G|-k+1$. Let $B$ be a basis of $\mathbb{F}_{q}^{k}$ and let $M_{G}^{\uparrow n}$ be a matrix whose rows are indexed by $\binom{G}{k-1}$, whose columns are indexed by ordered pairs $(U,A)$ where $U \in \binom{G}{n}$ and $A \in \binom{G \setminus U}{k-2}$, and whose $(C,(U,A))$-entry is
\begin{equation}
\label{entriesofMGr}
M_{G}^{\uparrow n}(C,(U,A)) = \begin{cases}
         \prod_{u \in U} \det(u,C)_{B} & \mbox{if $A \subset C \subset G \setminus U$} \\
         0 & \mbox{otherwise}.
         \end{cases}
\end{equation}
In \eqref{entriesofMGr}, $\det(u,C)_{B}$ denotes the determinant of the matrix whose first row is $u$ written with respect to the basis $B$ and whose last $k-1$ rows are the elements of $C$ written with respect to the basis $B$ in the order inherited from $G$.
\end{definition}

Although the matrices $M_{G}^{\uparrow n}$ may seem unfamiliar, we claim that they are related to inclusion matrices, which are well-studied in combinatorics. Recall that the inclusion matrix $I_{r}(a,b)$ has its rows indexed by $\binom{\{1, \ldots, r\}}{a}$, its columns indexed by $\binom{\{1, \ldots, r\}}{b}$, and $(A,B)$-entry
\begin{equation}
\label{entriesinclusion}
I_{r}(a,b)_{(A,B)} = \begin{cases}
                        1 & \text{if $B \subset A$} \\
                        0 & \text{otherwise}. \\
                        \end{cases}
\end{equation}
For example, when $n=0$, the matrix $M_{G}^{\uparrow 0}$ is the inclusion matrix $I_{|G|}(k-1,k-2)$. When $n > 0$, the matrix $M_{G}^{\uparrow n}$ is formed by gluing together matrices which are equivalent to inclusion matrices. For a fixed $U \in \binom{G}{n}$, let $D_{U}$ be a diagonal matrix whose rows and columns are indexed by $\binom{G \setminus U}{k-1}$ and whose $(C,C)$-entry is $\prod_{u \in U} \det(u,C)_{B}$. We then have that the submatrix $M_{G}^{\uparrow n}(U)$ of $M_{G}^{\uparrow n}$ whose rows are indexed by $\binom{G \setminus U}{k-1}$ and whose columns are indexed by ordered pairs $(U,A)$, where $A \in \binom{G \setminus U}{k-2}$, equals $D_{U}I_{|G \setminus U|}(k-1,k-2)$.

It is easy to see that linear equivalence of the arcs $G$ and $G'$ induces equivalence of the corresponding matrices $M_{G}^{\uparrow n}$ and $M_{G'}^{\uparrow n}$. More precisely, if $G, G' \subset \mathbb{F}_{q}^{k}$ are linearly equivalent arcs and $B$ and $B'$ are the bases of $\mathbb{F}_{q}^{k}$ used in the construction of the matrices $M_{G}^{\uparrow n}$ and $M_{G'}^{\uparrow n}$ respectively, then there exist invertible matrices $N_{1}$ and $N_{2}$ so that $M_{G}^{\uparrow n} = N_{1}M_{G'}^{\uparrow n}N_{2}$. 

Our main results relate algebraic properties of the matrix $M_{G}^{\uparrow n}$ to properties of the arc $G$. For example, our first main result says that if $G$ is an arc whose matrix $M_{G}^{\uparrow n}$ has full row rank, then $G$ cannot be extended to a larger arc of a specific size.

\begin{theorem} 
\label{thm:mainnew}
Let $G \subset \mathbb{F}_{q}^{k}$ be an arc and let $n \in \mathbb{N}$ be a natural number such that
\begin{equation}
\label{mainnewconstraints}
n+k-1 \leq |G| \leq \frac{q+2k-2+n}{2}. 
\end{equation}
If the matrix $M_{G}^{\uparrow n}$ has full row rank, then the arc $G$ cannot be extended to an arc of size $q+2k-1+n-|G|$. 
\end{theorem}
\noindent The left-hand and right-hand sides of \eqref{mainnewconstraints} respectively are required so that the matrix $M_{G}^{\uparrow n}$ exists and so that the arc $G$ has size strictly smaller than $q+2k-1+n-|G|$.

Suppose $0 \leq n \leq q-2k+4$ so that $2k-3+n \leq q+1$. Also, suppose we can show that for all arcs $G \subset \mathbb{F}_{q}^{k}$ of size $2k-3+n$, the matrix $M_{G}^{\uparrow n}$ has full row rank. If an arc of size $q+2$ exists in $\mathbb{F}_{q}^{k}$, then it would contain a subarc $G$ of size $2k-3+n$ that can be extended to an arc of size $q+2k-1+n-|G|$, which contradicts \tref{mainnew}. Consequently, \tref{mainnew} allows us to eliminate the existence of arcs of size $q+2$ in $\mathbb{F}_{q}^{k}$ by proving that for all arcs $G \subset \mathbb{F}_{q}^{k}$ of size $2k-3+n$, the matrix $M_{G}^{\uparrow n}$ has full row rank.

\begin{corollary}
\label{cor:maincor}
If $0 \leq n \leq q - 2k + 4$ and for every arc $G \subset \mathbb{F}_{q}^{k}$ of size $2k-3+n$, the matrix $M_{G}^{\uparrow n}$ has full row rank, then $g(k,q) = q+1$. 
\end{corollary}

Since the matrices $M_{G}^{\uparrow n}$ are related to inclusion matrices, knowing the ranks of inclusion matrices over $\mathbb{F}_{q}$ will be crucial to verifying the condition in \cref{maincor}.

\begin{theorem}[Frankl \cite{Franklprank}, Wilson \cite{Wilsonprank}]
\label{thm:FWprankformula}
For fixed integers $0 \leq b \leq a \leq r-b$ and a prime $p = \ch(\mathbb{F}_{q})$, we have
\begin{equation}
\label{prankformula}
\rank_{\mathbb{F}_{q}} \; I_{r}(a,b) = \sum_{\substack{0 \leq i \leq b \\ p \nmid \binom{a-i}{b-i}}} \binom{r}{i} - \binom{r}{i-1}.
\end{equation}
\end{theorem}

\noindent For example, when $n=0$ and $G \subset \mathbb{F}_{q}^{k}$ is an arc of size $2k-3$, the matrix $M_{G}^{\uparrow 0}$ is the inclusion matrix $I_{2k-3}(k-1,k-2)$. \tref{FWprankformula} thus implies the first assertion of \tref{nzero}.

\begin{theorem}
\label{thm:nzero}
If $G \subset \mathbb{F}_{q}^{k}$ is an arc of size $2k-3$, then the matrix $M_{G}^{\uparrow 0}$ has full row rank exactly when $k \leq p$. Hence $g(k,q) = q+1$ when $k \leq p$.
\end{theorem}

\noindent The second assertion of \tref{nzero} follows from \cref{maincor} when $q$ is not prime. If $q$ is prime, then \cref{maincor} implies that $g(k,q) = q+1$ when $k \leq (q+4)/2$ and hence the second assertion of \tref{nzero} follows from duality. The second assertion of \tref{nzero} was first proved by Ball \cite{MDSprime}. 

In \sref{rankP1}, we again use \tref{FWprankformula} to verify the condition in \cref{maincor} when $n=1$ and $k \leq 2p-2 \leq q$. 

\begin{theorem}
\label{thm:none}
If $k \leq 2p-2 \leq q$ and $G \subset \mathbb{F}_{q}^{k}$ is an arc of size $2k-2$, then the matrix $M_{G}^{\uparrow 1}$ has full row rank. Hence, if $q$ is not prime, then $g(k,q) = q+1$ when $k \leq 2p-2$.
\end{theorem}

\noindent The bound $k \leq 2p-2$ in the first assertion of \tref{none} cannot be improved because one can check using a computer that if $G \subset \mathbb{F}_{9}^{5}$ is a subarc of size $8$ of the normal rational curve $\cR_{5} \subset \mathbb{F}_{9}^{5}$, then the matrix $M_{G}^{\uparrow 1}$ does not have full row rank. The second assertion of \tref{none} follows from \cref{maincor} and was first proved by Ball and De Buele \cite{MDS2pminus2}. 

Recalling that $p = \ch(\mathbb{F}_{q})$, we conjecture that if $0 \leq n \leq q$ and
\begin{equation}
\label{rangek}
2 \leq k \leq \min \left \{p+n(p-2), \frac{q+4-n}{2} \right \},
\end{equation}
then the condition in \cref{maincor} holds. 

\begin{conjecture}
\label{conj:PErrank}
If $0 \leq n \leq q$, $k$ satisfies \eqref{rangek}, and $G \subset \mathbb{F}_{q}^{k}$ is an arc of size $2k-3+n$, then the matrix $M_{G}^{\uparrow n}$ has full row rank.
\end{conjecture}

\noindent Observe that \tref{nzero} and \tref{none} prove \conjref{PErrank} when $n=0$ and $n=1$. For larger values of $n$, we have computational evidence to support \conjref{PErrank}.

If \conjref{PErrank} is true then, by \cref{maincor}, $g(k,q)=q+1$ when \eqref{MDSbound} holds.

\begin{corollary}
\label{cor:goal}
If \conjref{PErrank} is true for any particular $n$ satisfying 
\begin{equation}
\label{rangen}
0 \leq n \leq \left| \frac{q-2p+4}{2p-3} \right|,
\end{equation}
then $g(k,q)=q+1$ when $k \leq p + n(p-2)$.  If \conjref{PErrank} is true, then $g(k,q) = q+1$ when \eqref{MDSbound} holds.
\end{corollary}

\subsubsection{Classification}
\label{sec:classificationresults}

The matrices $M_{G}^{\uparrow n}$ are also useful for determining when the normal rational curve $\cR_{k} \subset \mathbb{F}_{q}^{k}$ is the unique arc of size $q+1$ up to linear equivalence. The second main result of this article is that if $0 \leq n \leq q-2k$ and for any arc $G \subset \mathbb{F}_{q}^{k}$ of size $2k-2+n$, the matrix $M_{G}^{\uparrow n}$ contains a certain vector in its column space, then the normal rational curve is the unique arc of size $q+1$ up to linear equivalence. 

To state our theorem precisely, we define a matrix $H_{G}^{\uparrow n}$ that is equivalent to the matrix $M_{G}^{\uparrow n}$ so that the vector we require in the column space has a nice form. Recall that we have defined arcs to be ordered sets and that if $(X, <)$ is an ordered set then $A \subset X$ is smaller than $B \subset X$ in colex order if the largest element of the symmetric difference $A \triangle B$ lies in $B$. 

\begin{definition}
\label{def:ratiodiagonal}
Let $G \subset \mathbb{F}_{q}^{k}$ be an arc, let $0 \leq n \leq |G|-k+1$, and let $B$ be the basis of $\mathbb{F}_{q}^{k}$ fixed in \dref{M}. For each $C \in \binom{G}{k-1}$, let $L_{C} \in \binom{G \setminus C}{n}$ be the last $n$-subset of $\binom{G \setminus C}{n}$ in colex order. Let $J_{G}^{\uparrow n}$ be a diagonal matrix with rows and columns indexed by $\binom{G}{k-1}$ and $(C,C)$-entry
\begin{equation}
\label{entriesofJGn}
J_{G}^{\uparrow n}(C,C) = \prod_{y \in L_{C}} \det(y,C)_{B}^{-1}.
\end{equation}
Define the matrix $H_{G}^{\uparrow n} = J_{G}^{\uparrow n}M_{G}^{\uparrow n}$ and put the rows of the matrix $H_{G}^{\uparrow n}$ in colex order.
\end{definition}

Observe that the entries of the matrix $H_{G}^{\uparrow n}$ are independent of the basis $B$. We restate our second main result precisely using the matrices $H_{G}^{\uparrow n}$.

\begin{theorem}
\label{thm:classificationgeneral}
If $0 \leq n \leq q-2k$ and for every arc $G \subset \mathbb{F}_{q}^{k}$ of size $2k-2+n$, the column space of the matrix $H_{G}^{\uparrow n}$ contains a vector $v \in \mathbb{F}_{q}^{\binom{2k-2+n}{k-1}}$ such that $v_{i} = 1$ if $i \in \{1, \ldots, k \}$ and $v_{i} = 0$ otherwise, then the normal rational curve $\cR_{k}$ is the unique arc in $\mathbb{F}_{q}^{k}$ of size $q+1$ up to linear equivalence.
\end{theorem}

When $n=0$ and $G \subset \mathbb{F}_{q}^{k}$ is an arc of size $2k-2$, the matrix $H_{G}^{\uparrow 0}$ equals the inclusion matrix $I_{2k-2}(k-1,k-2)$, so we can easily verify that the column space of the matrix $H_{G}^{\uparrow 0}$ contains the required vector when $k \leq p = \ch(\mathbb{F}_{q})$.

\begin{theorem}
\label{thm:classificationzero}
If $k \leq p = \ch(\mathbb{F}_{q})$ and $G \subset \mathbb{F}_{q}^{k}$ is an arc of size $2k-2$, then the column space of the matrix $H_{G}^{\uparrow 0}$ contains a vector $v \in \mathbb{F}_{q}^{\binom{2k-2}{k-1}}$ such that $v_{i} = 1$ if $i \in \{1, \ldots, k \}$ and $v_{i} = 0$ otherwise. Hence, if $k \leq p$ and $k \neq (q+1)/2$, then the normal rational curve $\cR_{k}$ is the unique arc in $\mathbb{F}_{q}^{k}$ of size $q+1$ up to linear equivalence.
\end{theorem}

\noindent It is easy to see that the bound $k \leq p$ in the first assertion of \tref{classificationzero} cannot be improved. The second assertion of \tref{classificationzero} was first proved by Ball in \cite{MDSprime}, although the condition $k \neq (q+1)/2$ was missing there. 

We conjecture in \conjref{classify} that if $k \leq 2p-2 \leq q$ and $G \subset \mathbb{F}_{q}^{k}$ is an arc of size $2k$, then the column space of the matrix $H_{G}^{\uparrow 2}$ contains the required vector in \tref{classificationgeneral}. We have computational evidence to support \conjref{classify}, and we note that if \conjref{classify} is true, then the normal rational curve $\cR_{k}$ is the unique arc in $\mathbb{F}_{q}^{k}$ of size $q+1$ up to linear equivalence when $k \leq 2p-2 \leq q$.

\begin{conjecture}
\label{conj:classify}
When $k \leq 2p-2 \leq q$, for every arc $G \subset \mathbb{F}_{q}^{k}$ of size $2k$, the column space of the matrix $H_{G}^{\uparrow 2}$ contains a vector $v \in \mathbb{F}_{q}^{\binom{2k}{k-1}}$ such that $v_{i} = 1$ if $i \in \{1, \ldots, k \}$ and $v_{i} = 0$ otherwise.
\end{conjecture}

\subsubsection{Verifying \conjref{main} and Classifying Large Arcs Computationally}
\label{sec:computationalresults}

An important benefit of the conditions in \cref{maincor} and \tref{classificationgeneral} is that they can be checked with a computer. \cref{maincor} and \tref{classificationgeneral} may consequently be of use in verifying \conjref{main} and classifying large arcs computationally. For example, if one could classify arcs in $\mathbb{F}_{q}^{k}$ of size $2k-2$ up to linear equivalence, then one could test the rank of the matrix $M_{G}^{\uparrow 1}$ for a representative $G$ from each linear equivalence class. If the matrix $M_{G}^{\uparrow 1}$ has full row rank, then \cref{maincor} would rule out the possibility that any arc in the linear equivalence class of $G$ could be extended to an arc of size $q+2$. If the matrix $M_{G}^{\uparrow 1}$ does not have full row rank, then one could extend $G$ to an arc $H$ of size $2k-1$ and check if the matrix $M_{H}^{\uparrow 2}$ has full row rank. This should dramatically reduce the space of possible subarcs of arcs of size $q+2$. In the same way, \tref{classificationgeneral} can be used to check if the normal rational curve $\cR_{k}$ is the unique arc in $\mathbb{F}_{q}^{k}$ of size $q+1$ up to linear equivalence. These algorithms should be possible to implement because the question of classifying arcs up to linear equivalence has already been considered in \cite{Gordon} and \cite{Keri}.

\subsection{Important Remarks and Outline of Paper}
\label{sec:outline}

The results in this paper are joint work with Simeon Ball, but he has elected to write a separate exposition of some of these results in \cite{Ballextension}. A straightforward consequence of the proof of \tref{mainnew} is \tref{mainnewweightone}, which shows that that the conclusion of \tref{mainnew} holds if the matrix $M_{G}^{\uparrow n}$ satisfies the slightly weaker condition of having a vector of weight one in its column space. \tref{mainnewweightone} is the main result of \cite{Ballextension}.

\begin{theorem}
\label{thm:mainnewweightone}
Let $G \subset \mathbb{F}_{q}^{k}$ be an arc and let $n \in \mathbb{N}$ be a natural number such that
\begin{equation}
\label{mainnewconstraintsweightone}
n+k-1 \leq |G| \leq \frac{q+2k-2+n}{2}. 
\end{equation}
If the matrix $M_{G}^{\uparrow n}$ has a vector of weight one in its column space, then the arc $G$ cannot be extended to an arc of size $q+2k-1+n-|G|$. 
\end{theorem}

For the most interesting application of \tref{mainnew}, namely \cref{maincor}, we do not believe that \tref{mainnewweightone} offers any benefit over \tref{mainnew}. In other words, we believe that if $0 \leq n \leq q-2k+4$ and if for every arc $G \subset \mathbb{F}_{q}^{k}$ of size $2k-3+n$ the matrix $M_{G}^{\uparrow n}$ has a vector of weight one in its column space, then for every such arc $G$ the matrix $M_{G}^{\uparrow n}$ has full row rank. Indeed, the bound on $k$ in our stronger \conjref{PErrank} matches exactly the bound on $k$ in Ball's weaker Conjecture 1 in \cite{Ballextension}.

This paper builds on the methods initiated in \cite{MDSprime}, \cite{Ballbook}, and \cite{MDS2pminus2}. In order for this paper to be self-contained and correct, we repeat some proofs from \cite{MDSprime}, \cite{Ballbook}, and \cite{MDS2pminus2}, although we often give different expositions using matrices so that we may extend the results. For example, the proofs of \cite[Lemma 4.2]{MDSprime}, \cite[Lemma 7.20]{Ballbook}, and \cite[Lemma 3.1]{MDS2pminus2} are incorrect as written because the determinants in those results are evaluated with respect to many different bases yet treated as if they were evaluated with respect to the same basis. We fix this here in \cref{extension} and in \lref{mindlessr} and in \cite[Section 2]{Ballextension}. The proof of \cite[Theorem 1.8]{MDSprime} is also incorrect as written in the case $k = (q+1)/2$, and this is fixed here in \tref{classificationzero}.

Another important change in the proof approach of \cite{MDSprime}, \cite{Ballbook}, and \cite{MDS2pminus2} lies in the definition of certain parameters $\alpha_{A}$ in \lref{nullL}. In \cite[Chapter 7]{Ballbook}, the analogue of the parameter $\alpha_{A}$ in \lref{nullL} is referred to as $Q(A,F)$ and its definition is dependent on a smaller subarc of a larger arc. In \lref{nullL} and in \cite[Section 3]{Ballextension}, we define the parameters $\alpha_{A}$ so that they no longer depend on the smaller subarc and only depend on the larger arc. This change is crucial to the proof of \tref{mainnew}.   

The three main ingredients in the proofs of \tref{mainnew} and \tref{classificationgeneral} are duality, polynomial interpolation, and Segre's Lemma of Tangents. \sref{polyinterpolate} discusses the properties of polynomial interpolation that we use. \sref{tangent} explains the concept of tangent functions. In \sref{mainnewsketch}, we reduce our first main result \tref{mainnew} to \tref{similarr}. In \sref{similarrsketch}, we reduce \tref{similarr} to \lref{crux}. In \sref{segre}, we state and prove Segre's Lemma of Tangents and use it to prove \lref{crux}, thus completing the proofs of \tref{mainnew} and \tref{similarr}. In \sref{rankP1}, we prove \tref{none} and thus prove \conjref{PErrank} when $n=1$. In \sref{classification}, we prove \tref{classificationgeneral} and \tref{classificationzero}. 

\section{Polynomial Interpolation}
\label{sec:polyinterpolate}

That one can uniquely determine a polynomial $f \in \mathbb{F}[X]$ in one variable of degree at most $t$ over any field $\mathbb{F}$ from $t+1$ of its values is well-known. Similarly, one can recover a homogeneous polynomial in two variables $f \in \mathbb{F}(X,Y)$ of degree $t$ by knowing values of $f$ on the points of an arc $\{ (x_{i}, y_{i}) : i \in \{1, \ldots, t+1\}\}$ of size $t+1$ in $\mathbb{F}^{2}$.

Suppose $f(X,Y) = \sum_{i=0}^{t} c_{i}X^{i}Y^{t-i}$ is a homogeneous polynomial in two variables of degree $t$ and we know its values $f(x_{i}, y_{i})$ on the points of an arc $\{(x_{i}, y_{i}) : i \in \{1, \ldots, t+2\}\}$ of size $t+2$ in $\mathbb{F}^{2}$. Let $P \in M_{t+1, t+2}(\mathbb{F})$ be a matrix with $(i,j)$-entry $P(i,j) = x_{j}^{i-1}y_{j}^{t-i+1}$ and let $\vec{c} = [c_{0}, \ldots, c_{t}]$ and $\vec{z} = [f(x_{1}, y_{1}), \ldots, f(x_{t+2}, y_{t+2})]$. As $P$ has more columns than rows, its columns are linearly dependent. Hence, there is a solution $\vec{w} = [w_{1}, \ldots, w_{t+2}]^{T}$ to $P \vec{w} = \vec{0}$ and thus $\vec{z} \vec{w} = \vec{0}$ because $\vec{c} P = \vec{z}$. We now show in \tref{interpolationeqns} that a solution $\vec{w}$ to $P \vec{w} = \vec{0}$ and $\vec{z} \vec{w} = \vec{0}$ is given by 
\begin{equation}
\label{wsoln}
w_{i} = \prod_{\substack{j = 1 \\ j \neq i}}^{t+2} (x_{i}y_{j} - x_{j}y_{i})^{-1}, \; \; \; i \in \{1, \ldots, t+2\}.
\end{equation}
\tref{interpolationeqns} is a key ingredient in the proof of \tref{mainnew}.

\begin{theorem}
\label{thm:interpolationeqns}
Suppose $f(X,Y) \in \mathbb{F}[X,Y]$ is a homogeneous polynomial in two variables of degree $t$ and $\{(x_{i}, y_{i}) : i \in \{1, \ldots, t+2\}\}$ is an arc of size $t+2$ in $\mathbb{F}^{2}$. We then have
\begin{equation}
\label{interpolationeqn}
\sum_{i=1}^{t+2} f(x_{i}, y_{i}) \prod_{\substack{j = 1 \\ j \neq i}}^{t+2} (x_{i}y_{j} - x_{j}y_{i})^{-1}  = 0.
\end{equation}
\end{theorem}

\noindent \textbf{Proof.} Using the definitions of $P$, $\vec{w}$, and $\vec{z}$ from the preceding paragraph, let $B$ be a square matrix whose columns are the first $t+1$ columns of the matrix $P$. Let $\vec{b}$ be the last column of the matrix $P$. Note that a solution $\vec{r} = [r_{1}, \ldots, r_{t+1}]^{T}$ to $B \vec{r} = \vec{b}$ gives a solution $\vec{w}$ to $P \vec{w} = \vec{0}$ with $w_{i} = r_{i}$ for $i \in \{1, \ldots, t+1\}$ and $w_{t+2} = -1$.

Since $\{(x_{i}, y_{i}) : i \in \{1, \ldots, t+2\}\}$ is an arc of size $t+2$ in $\mathbb{F}^{2}$, we may assume that $y_{1}, \ldots, y_{t+1}$ are nonzero. Hence the matrix $B$ is nonsingular, so by Cramer's Rule, a solution $\vec{r}$ to $B \vec{r} = \vec{b}$ is given by $r_{i} = \det(B_{i})/\det(B)$ where $B_{i}$ is the matrix formed by replacing the $i^{\mathrm{th}}$ column of $B$ with $\vec{b}$. Using the formula for the determinant of a Vandermonde matrix, 
\begin{align*}
\det(B) &= (y_{1} \cdots y_{t+1})^{t} \prod_{1 \leq l < m \leq t+1} \left( \frac{x_{m}}{y_{m}} - \frac{x_{l}}{y_{l}} \right) = \prod_{1 \leq l < m \leq t+1} (x_{m}y_{l} - x_{l}y_{m}), \\
\det(B_{i}) &= \prod_{\substack{1 \leq l < m \leq t+1 \\ l \neq i, m \neq i}} (x_{m}y_{l} - x_{l}y_{m}) \prod_{1 \leq l < i} (x_{t+2}y_{l} - x_{l}y_{t+2}) \prod_{i < m \leq t+1} (x_{m}y_{t+2} - x_{t+2}y_{m}).  
\end{align*}
Hence, after a little algebraic manipulation,
\begin{equation}
\label{applyCramer}
r_{i} = \frac{\det(B_{i})}{\det(B)} = -\frac{ \prod_{1 \leq l \leq t+1} (x_{t+2}y_{l} - x_{l}y_{t+2})}{ \prod_{\substack{1 \leq l \leq t+2 \\ l \neq i}} (x_{i}y_{l} - x_{l}y_{i})}.  
\end{equation}
Multiplying the corresponding solution $\vec{w}$ to $P \vec{w} = \vec{0}$ by $- \prod_{1 \leq l \leq t+1} (x_{t+2}y_{l} - x_{l}y_{t+2})^{-1}$ yields the solution $\vec{w}$ to $\vec{z} \vec{w} = \vec{0}$ given by \eqref{wsoln}. \qed

\section{Tangent Functions}
\label{sec:tangent}

Let $S \subset \mathbb{F}_{q}^{k}$ be an arc. Given a subset $A \subset S$ of size $k-2$, we will define the tangent function at $A$, denoted $f_{A,S}: \mathbb{F}_{q}^{k} \rightarrow \mathbb{F}_{q}$, which can be viewed as a homogeneous polynomial in two variables with respect to certain bases of $\mathbb{F}_{q}^{k}$. We will then apply \tref{interpolationeqns} to the tangent functions $f_{A,S}$ for various $A \subset S$ to prove \tref{mainnew}. 

To define the tangent function at $A$, we first count in \lref{degreetangent} the number of $(k-1)$-dimensional subspaces of $\mathbb{F}_{q}^{k}$ that intersect $S$ precisely in $A$. 

\begin{lemma}[Ball \cite{MDSprime, Ballbook}] 
\label{lem:degreetangent}
Let $S \subset \mathbb{F}_{q}^{k}$ be an arc and let $A \subset S$ have size $k-2$. Let $H_{A}^{1}, \ldots, H_{A}^{t}$ be the $(k-1)$-dimensional subspaces of $\mathbb{F}_{q}^{k}$ whose intersection with $S$ is $A$. We have 
\begin{equation}
\label{defoft}
t := q+k-1 - |S|. 
\end{equation} 
\end{lemma}

\noindent \textbf{Proof.} Since $A$ is a linearly independent set of size $k-2$, the number of  $(k-1)$-dimensional subspaces of $\mathbb{F}_{q}^{k}$ that contain $A$ is $q+1$. Since $S$ is an arc, a $(k-1)$-dimensional subspace of $\mathbb{F}_{q}^{k}$ that contains $A$ can contain at most one other vector of $S \setminus A$. \qed

\vspace{0.25cm}

Given an arc $S \subset \mathbb{F}_{q}^{k}$ and a subset $A \subset S$ of size $k-2$, we now define the tangent function at $A$.

\begin{definition}[Ball \cite{MDSprime, Ballbook}]
\label{def:tangentfunction}
Let $S \subset \mathbb{F}_{q}^{k}$ be an arc and let $A \subset S$ have size $k-2$. Let $H_{A}^{1}, \ldots, H_{A}^{t}$ be the $(k-1)$-dimensional subspaces of $\mathbb{F}_{q}^{k}$ defined in \lref{degreetangent}, where $t$ is given by \eqref{defoft}. Let $\beta_{A}^{i}: \mathbb{F}_{q}^{k} \rightarrow \mathbb{F}_{q}$ be a linear functional whose kernel is $H_{A}^{i}$. We define the tangent function at $A$, denoted $f_{A,S} : \mathbb{F}_{q}^{k} \rightarrow \mathbb{F}_{q}$, by
\begin{equation}
\label{tangentfunctiondef}
f_{A,S}(x) = \prod_{i=1}^{t} \beta_{A}^{i}(x).
\end{equation}
\end{definition}

\noindent Observe that $f_{A,S}(x) = 0$ precisely when $x \in \bigcup_{i=1}^{t} H_{A}^{i}$ and that $f_{A,S}$ is defined up to a scalar factor.

\vspace{0.25cm}

\textbf{Notation:} Recall that an arc $S \subset \mathbb{F}_{q}^{k}$ is ordered. If $R_{1}, \ldots, R_{l}$ are subsets of $S$ we use $(R_{1}, \ldots, R_{l})$ to mean write the vectors in $R_{1}$ in order first, and then the vectors in $R_{2}$ etc. When $R_{i}$ is a singleton set, we simply write the vector. For example, if $x, y \in S \setminus A$ and $B$ is a basis of $\mathbb{F}_{q}^{k}$, we write $\det(x,y,A)_{B}$ for the determinant of the matrix whose rows are the vectors $x$, $y$, and the elements of $A$ in order written with respect to $B$.  

\vspace{0.25cm}

Let $S \subset \mathbb{F}_{q}^{k}$ be an arc and let $A \subset S$ have size $k-2$. Let $B = (b_{1}, b_{2}, A)$ be a basis of $\mathbb{F}_{q}^{k}$. Also suppose that $T \subset S \setminus A$ is a subset of size $t+2$, where $t$ is defined by \eqref{defoft}. In \lref{homogtwovar}, we show how to use \tref{interpolationeqns} to obtain an equation for a pair $(A,T)$ where $A \in \binom{S}{k-2}$ and $T \in \binom{S \setminus A}{t+2}$. Such equations are crucial to the proof of \tref{mainnew}.

\begin{lemma}[Ball \cite{MDSprime, Ballbook}]
\label{lem:homogtwovar}
Let $S \subset \mathbb{F}_{q}^{k}$ be an arc and let $A \subset S$ have size $k-2$. Let $B = (b_{1}, b_{2}, A)$ be a basis of $\mathbb{F}_{q}^{k}$. If $T \subset S \setminus A$ has size $t+2$, where $t$ is given by \eqref{defoft}, then 
\begin{equation}
\label{howtoapplyinterpolation}
\sum_{x \in T} f_{A,S}(x) \prod_{y \in T \setminus \{x\}} \det(x,y,A)^{-1}_{B} = 0.
\end{equation}
\end{lemma}

\noindent \textbf{Proof.} With respect to the basis $B$, the linear functional $\beta_{A}^{i}$ in \eqref{tangentfunctiondef} is linear in just the first two coordinates since its kernel contains $A$. Hence, the tangent function $f_{A,S}$ is a homogeneous polynomial in two variables of degree $t$, where $t$ is given by \eqref{defoft}. Since $S$ is an arc, when we write the vectors in $T$ in terms of the basis $B$, their first two coordinates form an arc of size $t+2$ in $\mathbb{F}_{q}^{2}$. Hence, we can apply \tref{interpolationeqns} to $f_{A,S}$ and $T$, and note that with respect to $B$, we have $\det(x,y,A)_{B} = x_{1}y_{2} - y_{1}x_{2}$. \qed

We will show in \lref{mindlessr} that the product of determinants in \eqref{howtoapplyinterpolation} is related to the product of determinants in the entries of the matrix $M_{G}^{\uparrow n}$ in \eqref{entriesofMGr}. Hence \tref{interpolationeqns} and \lref{homogtwovar} explain how the entries of the matrix $M_{G}^{\uparrow n}$ arise.

\section{Proof that \tref{similarr} Implies \tref{mainnew}}
\label{sec:mainnewsketch}

Let $S \subset \mathbb{F}_{q}^{k}$ be an arc and let $G \subset S$ have size $t+k+n$, where $t$ is defined by \eqref{defoft} and $n \geq 0$. For each $U \in \binom{G}{n}$, the system of equations obtained from applying \lref{homogtwovar} to pairs $(A,T)$ where $A \in \binom{G \setminus U}{k-2}$ and $T = G \setminus (A \cup U)$ can be readily analyzed and is a key ingredient in the proof of \tref{mainnew}. 

\begin{corollary}
\label{cor:extension}
Let $S \subset \mathbb{F}_{q}^{k}$ be an arc and let $G \subset S$ have size $t+k+n$, where $t$ is defined by \eqref{defoft} and $n \geq 0$. For $A \in \binom{G}{k-2}$, define a basis $B(A) = (b_{1}, b_{2}, A)$ of $\mathbb{F}_{q}^{k}$. If $P_{G}^{\uparrow n}$ is a matrix whose rows are indexed by $\binom{G}{k-1}$, whose columns are indexed by ordered pairs $(U,A)$ where $U \in \binom{G}{n}$ and $A \in \binom{G \setminus U}{k-2}$, and whose $(C,(U,A))$-entry is
\begin{equation*}
P_{G}^{\uparrow n}(C,(U,A)) = \begin{cases}
         f_{A,S}(C \setminus A) \prod_{y \in G \setminus (C \cup U)} \det(C \setminus A,y,A)^{-1}_{B(A)} & \mbox{if $A \subset C \subset G \setminus U$} \\
         0 & \mbox{otherwise},
         \end{cases}
\end{equation*}
then $\vec{1} P_{G}^{\uparrow n} = \vec{0}$.
\end{corollary}

\noindent \textbf{Proof.} For $U \in \binom{G}{n}$ and $A \in \binom{G \setminus U}{k-2}$, let $T = G \setminus (A \cup U)$. Since $|G|=t+k+n$, we have that $T \in \binom{S \setminus A}{t+2}$. Applying \lref{homogtwovar}, we have that 
\begin{equation}
\label{tangentinterpolationeqnr}
\sum_{x \in G \setminus (A \cup U)} f_{A,S}(x) \prod_{y \in (G \setminus (A \cup U)) \setminus \{x\}} \det(x,y,A)^{-1}_{B(A)} = 0.
\end{equation}
Let us rewrite \eqref{tangentinterpolationeqnr} so that it will be easier to express the system of equations given by \eqref{tangentinterpolationeqnr} in matrix form. For any fixed $U \in \binom{G}{n}$ and $A \in \binom{G \setminus U}{k-2}$, we have
\begin{equation}
\label{tangentinterpolaterewriter}
\sum_{A \subset C \in \binom{G \setminus U}{k-1}} f_{A,S}(C \setminus A) \prod_{y \in G \setminus (C \cup U)} \det(C \setminus A,y,A)^{-1}_{B(A)} = 0.
\end{equation}
Consequently, letting $P_{G}^{\uparrow n}$ be the matrix defined in \cref{extension}, we see that we can write the system of equations given by \eqref{tangentinterpolaterewriter} in matrix form as $\vec{1} P_{G}^{\uparrow n} = \vec{0}$. \qed

The equation $\vec{1} P_{G}^{\uparrow n} = \vec{0}$ contains a wealth of information about the arc $S \subset \mathbb{F}_{q}^{k}$ and is crucial to the proof of \tref{mainnew}. At the moment, the matrix $P_{G}^{\uparrow n}$ defined in \cref{extension} may seem ugly and difficult to analyze, but we claim that $P_{G}^{\uparrow n}$ is equivalent to the much simpler matrix $M_{G}^{\uparrow n}$ defined in \eqref{entriesofMGr}, which depends only on the arc $G$.

\begin{theorem}
\label{thm:similarr}
Let $S \subset \mathbb{F}_{q}^{k}$ be an arc and let $G \subset S$ have size $t+k+n$, where $t$ is defined by \eqref{defoft} and $n \geq 0$. If $P_{G}^{\uparrow n}$ is the matrix defined in \cref{extension}, then there exist invertible diagonal matrices $D_{1}$ and $D_{2}$ so that $D_{1} P_{G}^{\uparrow n} D_{2} = M_{G}^{\uparrow n}$, where $M_{G}^{\uparrow n}$ is defined by \eqref{entriesofMGr}.
\end{theorem}

We now reduce \tref{mainnew} and \tref{mainnewweightone} to \tref{similarr}. 

\vspace{0.25cm}

\noindent \textbf{Proof of \tref{mainnew} and \tref{mainnewweightone}}. We prove the contrapositive: namely that if $G \subset \mathbb{F}_{q}^{k}$ can be extended to an arc $S \subset \mathbb{F}_{q}^{k}$ of size $q+2k-1+n-|G|$, then the matrix $M_{G}^{\uparrow n}$ cannot have full row rank or a vector of weight one in its column space. First we show the arc $G$ satisfies the hypotheses of \cref{extension}. As $S$ has size $q+2k-1+n-|G|$, the arc $G$ has size $t+k+n$, where $t$ is defined by \eqref{defoft}. By \cref{extension}, we have $\vec{1} P_{G}^{\uparrow n} = \vec{0}$ and so by \tref{similarr}, we have $\vec{0} = (\vec{1}D_{1}^{-1})M_{G}^{\uparrow n}$. Since $D_{1}$ is an invertible matrix, all entries of $\vec{1} D_{1}^{-1}$ are nonzero. Hence, the matrix $M_{G}^{\uparrow n}$ cannot have full row rank or a vector of weight one in its column space. \qed

\section{Proof that \lref{crux} Implies \tref{similarr}}
\label{sec:similarrsketch}

To prove \tref{similarr}, we first define matrices $Q_{G}^{\uparrow n}$, $R_{G}^{\uparrow n}$ and $I_{G}^{\uparrow n}$ given a subset $G \subseteq S$ where $S \subset \mathbb{F}_{q}^{k}$ is an arc. The matrix $I_{G}^{\uparrow n}$ is a signed inclusion matrix and to define the signing we need the following notation.

\begin{definition}
\label{def:tau}
Let $X$ be an ordered set, let $A$ be an ordered subset of $X$, and let $C$ be an ordered subset of $X$ that contains $A$ and has size $|A|+1$. We define $\tau(A,C)$ to be the minimum number of transpositions needed to order $(A, C \setminus A)$ as $C$.
\end{definition}

\begin{definition}
\label{def:Q}
Let $G \subseteq S \subset \mathbb{F}_{q}^{k}$, where $S$ is an arc, and let $0 \leq n \leq |G|-k+1$. Let $Q_{G}^{\uparrow n}$, $R_{G}^{\uparrow n}$, and $I_{G}^{\uparrow n}$ respectively be matrices whose rows are indexed by $\binom{G}{k-1}$, whose columns are indexed by ordered pairs $(U,A)$ where $U \in \binom{G}{n}$ and $A \in \binom{G \setminus U}{k-2}$, and whose $(C,(U,A))$-entries respectively are
\begin{equation}
\label{entriesofQr}
Q_{G}^{\uparrow n}(C,(U,A)) = \begin{cases}
         f_{A,S}(C \setminus A) & \mbox{if $A \subset C \subset G \setminus U$}  \\
         0 & \mbox{otherwise},
         \end{cases} 
\end{equation}
\begin{equation}
\label{entriesofRr} 
R_{G}^{\uparrow n}(C,(U,A)) = \begin{cases}
         \prod_{y \in G \setminus (C \cup U)} \det(C \setminus A,y,A)^{-1}_{B(A)} & \mbox{if $A \subset C \subset G \setminus U$} \\
         0 & \mbox{otherwise},
         \end{cases} 
\end{equation}
\begin{equation}
\label{entriesinclusionliker}
I_{G}^{\uparrow n}(C,(U,A)) = \begin{cases}
                        (-1)^{\tau(A,C)(t+1)} & \text{if $A \subset C \subset G \setminus U$} \\
                        0 & \text{otherwise}, \\
                        \end{cases}
\end{equation}
where $f_{A,S}$ is defined by \dref{tangentfunction}, $\tau(A,C)$ is defined by \dref{tau}, and $t$ is defined by \eqref{defoft}. 
\end{definition}

We will prove in \sref{segre} that if $S \subset \mathbb{F}_{q}^{k}$ is an arc, then the matrix $Q_{S}^{\uparrow 0}$ defined in \eqref{entriesofQr} is equivalent to the matrix ${I}_{S}^{\uparrow 0}$ defined in \eqref{entriesinclusionliker}.

\begin{lemma}
\label{lem:crux}
Let $S \subset \mathbb{F}_{q}^{k}$ be an arc. If $Q_{S}^{\uparrow 0}$ is the matrix defined in \eqref{entriesofQr} and ${I}_{S}^{\uparrow 0}$ is the matrix defined in \eqref{entriesinclusionliker}, then there exist invertible diagonal matrices $E_{1}$ and $E_{2}$ such that $E_{1}Q_{S}^{\uparrow 0}E_{2} = {I}_{S}^{\uparrow 0}$.
\end{lemma}

We now use \lref{crux} to prove that if $S \subset \mathbb{F}_{q}^{k}$ is an arc and $G \subset S$ satisfies the constraints of \cref{extension} then the matrix $Q_{G}^{\uparrow n}$ defined in \eqref{entriesofQr} is equivalent to the matrix $I_{G}^{\uparrow n}$ defined in \eqref{entriesinclusionliker}.

\begin{lemma}
\label{lem:cruxr}
Let $G \subset S \subset \mathbb{F}_{q}^{k}$ where $S$ is an arc. If $Q_{G}^{\uparrow n}$ is the matrix defined in \eqref{entriesofQr}, then there exist invertible diagonal matrices $F_{1}$ and $F_{2}$ such that $F_{1}Q_{G}^{\uparrow n}F_{2} = I_{G}^{\uparrow n}$.
\end{lemma}

\noindent \textbf{Proof.} Recall that, by \lref{crux}, there exist invertible diagonal matrices $E_{1}$ and $E_{2}$ such that $E_{1}Q_{S}^{\uparrow 0}E_{2} = I_{S}^{\uparrow 0}$. Let $F_{1}$ be the submatrix of $E_{1}$ whose rows and columns are indexed by $\binom{G}{k-1}$. Let $F_{2}$ be the submatrix of $E_{2}$ whose rows and columns are indexed by ordered pairs $(U,A)$, where $U \in \binom{G}{n}$ and $A \in \binom{G \setminus U}{k-2}$. As the entries of $Q_{G}^{\uparrow n}$ and $I_{G}^{\uparrow n}$ don't depend on $U \in \binom{G}{n}$, we have $F_{1}Q_{G}^{\uparrow n}F_{2} = I_{G}^{\uparrow n}$. \qed

We now prove that if $S \subset \mathbb{F}_{q}^{k}$ is an arc and $G \subset S$ satisfies the constraints of \cref{extension} then the matrix $R_{G}^{\uparrow n}$ defined in \eqref{entriesofRr} is equivalent to the Hadamard product $I_{G}^{\uparrow n} \circ M_{G}^{\uparrow n}$.

\begin{lemma}
\label{lem:mindlessr}
Let $S \subset \mathbb{F}_{q}^{k}$ be an arc and let $G \subset S$ satisfy the hypotheses of \cref{extension}. If $R_{G}^{\uparrow n}$ is the matrix defined in \eqref{entriesofRr}, then there exist invertible diagonal matrices $F_{3}$ and $F_{4}$ such that $F_{3}R_{G}^{\uparrow n}F_{4} = I_{G}^{\uparrow n} \circ M_{G}^{\uparrow n}$.
\end{lemma}

\noindent \textbf{Proof.} Since $|G|=t+k+n$, for $C \in \binom{G}{k-1}$ and $U \in \binom{G}{n}$, we have $|G \setminus (C \cup U)| = t+1$. Hence, for $A \subset C$, we have
\begin{equation}
\label{movetoendr}
(-1)^{(t+1)(k-1)} \prod_{y \in G \setminus (C \cup U)} \det(C \setminus A, y, A)^{-1}_{B(A)} = \prod_{y \in G \setminus (C \cup U)} \det(y,A, C \setminus A)^{-1}_{B(A)} 
\end{equation}
because there are $k-1$ transpositions needed to make $C \setminus A$ the last row. 

Recall that for each $A \in \binom{G}{k-2}$, we defined a basis $B(A) = (b_{1}, b_{2}, A)$ of $\mathbb{F}_{q}^{k}$. Now let $B$ be the basis of $\mathbb{F}_{q}^{k}$ fixed in \dref{M} and let $M(B(A),B)$ be the change-of-basis matrix from $B(A)$ to $B$. Observe that 
\begin{equation}
\label{changeofbasis}
\det(y,A, C \setminus A)^{-1}_{B(A)} \det(M(B(A),B))^{-1} = \det(y, A, C \setminus A)^{-1}_{B}.
\end{equation}
Define $F_{4}$ to be a diagonal matrix with rows and columns indexed by ordered pairs $(U,A)$ where $U \in \binom{G}{n}$ and $A \in \binom{G \setminus U}{k-2}$, and $((U,A),(U,A))$-entry
\begin{equation}
\label{entriesofF4}
F_{4}((U,A),(U,A)) = (-1)^{(t+1)(k-1)} \det(M(B(A),B))^{-(t+1)}.
\end{equation}
By \eqref{movetoendr} and \eqref{changeofbasis}, the $(C,(U,A))$-entry of $R_{G}^{\uparrow n}F_{4}$ is 
\begin{equation}
\label{entriesofRF4}
R_{G}^{\uparrow n}F_{4}(C,(U,A)) = \begin{cases}
                      \prod_{y \in G \setminus (C \cup U)} \det(y, A, C \setminus A)^{-1}_{B} & \mbox{if $A \subset C \subset G \setminus U$} \\
                      0 & \mbox{otherwise}.
                      \end{cases}
\end{equation}

Observe that 
\begin{equation}
\label{detinclusionliker}
\prod_{y \in G \setminus (C \cup U)} \det(y, A, C \setminus A)^{-1}_{B} = (-1)^{\tau(A,C)(t+1)} \prod_{y \in G \setminus (C \cup U)} \det(y,C)^{-1}_{B}
\end{equation}
because moving $C \setminus A$ from the end to its proper place in the ordering of $C$ requires $\tau(A,C)$ transpositions for each of the $t+1$ determinants in the product.

Hence, defining $F_{3}$ to be a diagonal matrix with rows and columns indexed by $\binom{G}{k-1}$ and $(C,C)$-entry
\begin{equation}
\label{entriesofF3}
F_{3}(C,C) = \prod_{y \in G \setminus C} \det(y,C)_{B},
\end{equation}
we see that $F_{3}R_{G}^{\uparrow n}F_{4} = I_{G}^{\uparrow n} \circ M_{G}^{\uparrow n}$. \qed

Finally we reduce \tref{similarr} to \lref{crux}.

\vspace{0.25cm}

\noindent \textbf{Proof of \tref{similarr}} We express the matrix $P_{G}^{\uparrow n}$ defined in \cref{extension} as the Hadamard product $P_{G}^{\uparrow n} = Q_{G}^{\uparrow n} \circ R_{G}^{\uparrow n}$, where the matrices $Q_{G}^{\uparrow n}$ and $R_{G}^{\uparrow n}$ are defined in \eqref{entriesofQr} and \eqref{entriesofRr} respectively.

Using \lref{crux}, we show in \lref{cruxr} that there exist invertible diagonal matrices $F_{1}$ and $F_{2}$ such that $F_{1}Q_{G}^{\uparrow n}F_{2} = I_{G}^{\uparrow n}$. By \lref{mindlessr}, there exist invertible diagonal matrices $F_{3}$ and $F_{4}$ such that $F_{3}R_{G}^{\uparrow n}F_{4} = I_{G}^{\uparrow n} \circ M_{G}^{\uparrow n}$. Setting $D_{1} = F_{1}F_{3}$ and $D_{2} = F_{2}F_{4}$, \tref{similarr} follows. \qed

\section{Proof of \lref{crux}}
\label{sec:segre}

In this section, we prove \lref{crux} and hence complete the proofs of \tref{mainnew} and \tref{similarr}. In \lref{nullL}, we show that \lref{crux} holds if we can find a vector in the nullspace of a certain matrix $L$ all of whose coordinates are nonzero. In \lref{segrelemmaoftangents}, we state and prove Segre's Lemma of Tangents, which we use in \lref{rankL} to show that the matrix $L$ does not have full column rank. Consequently, $L$ has nonzero vectors in its nullspace and we prove \lref{crux} by showing that any nonzero vector in the nullspace of $L$ must have all coordinates nonzero.

Recall that an arc $S \subset \mathbb{F}_{q}^{k}$ is ordered and that if $(X, <)$ is an ordered set then $A \subset X$ is smaller than $B \subset X$ in lex order if the smallest element of the symmetric difference $A \triangle B$ lies in $A$. 

\begin{lemma}
\label{lem:nullL}
Let $S \subset \mathbb{F}_{q}^{k}$ be an arc. Let $L$ be a matrix whose columns are indexed by $\binom{S}{k-2}$ and whose rows are indexed by ordered pairs $(A,A')$ where $A,A' \in \binom{S}{k-2}$, $A \cup A' \in \binom{S}{k-1}$, and $A < A'$ in lex order. Let the $((A,A'),A'')$ entry of $L$ be 
\begin{equation}
\label{entriesofL}
L((A,A')),A'') = \begin{cases}
                (-1)^{\tau(A, A \cup A')(t+1)} f_{A,S}(A' \setminus A) & \mbox{if $A'' = A$} \\
                (-1)^{\tau(A', A \cup A')(t+1)+1} f_{A',S}(A \setminus A') & \mbox{if $A'' = A'$} \\
                0 & \mbox{otherwise},
                \end{cases}
\end{equation}
where $t$ is defined by \eqref{defoft}. If there exists a vector $\vec{\alpha} \in \mathbb{F}_{q}^{\binom{|S|}{k-2}}$ in the nullspace of $L$ all of whose coordinates are nonzero, then \lref{crux} holds.
\end{lemma}

\noindent \textbf{Proof.} We write the coordinates of $\vec{\alpha}$ as $\alpha_{A}$ where $A \in \binom{S}{k-2}$. Since $\vec{\alpha}$ is in the nullspace of $L$, if $Q_{S}^{\uparrow 0}$ is the matrix defined in \eqref{entriesofQr}, $C \in \binom{S}{k-1}$, $A, A' \in \binom{C}{k-2}$, and $t$ is defined by \eqref{defoft}, then
\begin{equation}
\label{rowentriesequal}
(-1)^{\tau(A,C)(t+1)} \alpha_{A} Q_{S}^{\uparrow 0}(C,A) = (-1)^{\tau(A',C)(t+1)} \alpha_{A'}Q_{S}^{\uparrow 0}(C,A').
\end{equation}
Define $E_{2}$ to be a diagonal matrix with rows and columns indexed by $\binom{S}{k-2}$ and $(A,A)$ entry $E_{2}(A,A) = \alpha_{A}$. Since the coordinates of $\vec{\alpha}$ are nonzero and $Q_{S}^{\uparrow 0}(C,A) \neq 0$ when $A \subset C$, there exist nonzero constants $\alpha_{C} \in \mathbb{F}_{q}$ for $C \in \binom{S}{k-1}$ such that the $(C,A)$ entry of $Q_{S}^{\uparrow 0} E_{2}$ is 
\begin{equation}
\label{entriesofQE2}
Q_{S}^{\uparrow 0}E_{2}(C,A) = \begin{cases}
              (-1)^{\tau(A,C)(t+1)} \alpha_{C} & \mbox{if $A \subset C$} \\
              0 & \mbox{otherwise},
              \end{cases}
\end{equation}
by \eqref{rowentriesequal}. Consequently, defining $E_{1}$ to be a diagonal matrix with rows and columns indexed by $\binom{S}{k-1}$ and $(C,C)$-entry $E_{1}(C,C) = \alpha_{C}^{-1}$, we see that \lref{crux} holds. \qed

To prove the existence of a vector $\vec{\alpha} \in \mathbb{F}_{q}^{\binom{|S|}{k-2}}$ satisfying the hypotheses of \lref{nullL}, we first show in \lref{rankL} that the matrix $L$ defined in \eqref{entriesofL} does not have full column rank over $\mathbb{F}_{q}$. For this, we need \lref{segrelemmaoftangents}, which is called Segre's Lemma of Tangents and gives a relationship between values of different tangent functions. 

\begin{lemma}[Ball \cite{MDSprime, Ballbook}]
\label{lem:segrelemmaoftangents}
Let $S \subset \mathbb{F}_{q}^{k}$ be an arc and let $t$ be defined by \eqref{defoft}. For a subset $D \subset S$ of size $k-3$ and a subset $\{u,v,w\} \in S \setminus D$, we have
\begin{equation}
\label{lemmaoftangents}
f_{D \cup \{u\},S}(v)f_{D \cup \{v\},S}(w) f_{D \cup \{w\},S}(u) = (-1)^{t+1} f_{D \cup \{u\},S}(w) f_{D \cup \{v\},S}(u) f_{D \cup \{w\},S}(v).
\end{equation}
\end{lemma}

\noindent \textbf{Proof.} Observe that $B = (u,v,w,D)$ is a basis of $\mathbb{F}_{q}^{k}$ because $S$ is an arc. For $x \in \mathbb{F}_{q}^{k}$, let $x=(x_{1}, \ldots, x_{k})$ be the coordinates of $x$ with respect to $B$. By \eqref{tangentfunctiondef},
\begin{equation}
\label{segretangent}
f_{D \cup \{w\},S}(x) = \prod_{i=1}^{t} (\beta_{D \cup \{w\}}^{i}(u) x_{1} + \beta_{D \cup \{w\}}^{i}(v) x_{2}).  
\end{equation}

Our first goal is to show that
\begin{equation}
\label{goal1}
\left \{-\frac{ \beta_{D \cup \{w\}}^{i}(u)}{\beta_{D \cup \{w\}}^{i}(v)} : i \in \{1, \ldots, t\} \right \} \cup \left \{ \frac{x_{2}}{x_{1}} : x \in S \setminus B \right \} = \mathbb{F}_{q} \setminus \{0\}.
\end{equation}
To accomplish this, observe that the first set on the left hand side of \eqref{goal1} contains $t$ nonzero elements of $\mathbb{F}_{q}$ because for $i \in \{1, \ldots, t\}$, the $(k-1)$-dimensional subspaces $H_{D \cup \{w\}}^{i}$ defined in \lref{degreetangent} are all distinct and
intersect $S$ only in $D \cup \{w\}$. Now observe that the second set on the left hand side of \eqref{goal1} is disjoint from the first set and contains $|S|-k$ nonzero elements of $\mathbb{F}_{q}$ because $S$ is an arc and because for $i \in \{1, \ldots, t\}$, the $(k-1)$-dimensional subspaces $H_{D \cup \{w\}}^{i}$ defined in \lref{degreetangent} intersect $S$ only in $D \cup \{w\}$.  Since $t + |S| - k = q-1$, \eqref{goal1} is established.

Since the product of the nonzero elements of a finite field $\mathbb{F}_{q}$ equals $-1$, \eqref{goal1} implies 
\begin{equation}
\label{minusone}
\prod_{i=1}^{t} \left(-\frac{\beta_{D \cup \{w\}}^{i}(u)}{\beta_{D \cup \{w\}}^{i}(v)} \right) \prod_{x \in S \setminus B} \frac{x_{2}}{x_{1}} = -1.
\end{equation}
By \eqref{segretangent}, we can rewrite \eqref{minusone} as
\begin{equation}
\label{replace}
f_{D \cup \{w\},S}(u) \prod_{x \in S \setminus B} x_{2} = (-1)^{t+1} f_{D \cup \{w\},S}(v) \prod_{x \in S \setminus B} x_{1}. 
\end{equation}
Repeating the argument above with the $(k-2)$-subsets $D \cup \{u\}$ and $D \cup \{v\}$, we have  
\begin{align}
f_{D \cup \{u\},S}(v) \prod_{x \in S \setminus B} x_{3} &= (-1)^{t+1} f_{D \cup \{u\},S}(w) \prod_{x \in S \setminus B} x_{2} \label{rinseandrepeat1} \\
f_{D \cup \{v\},S}(w) \prod_{x \in S \setminus B} x_{1} &= (-1)^{t+1} f_{D \cup \{v\},S}(u) \prod_{x \in S \setminus B} x_{3} \label{rinseandrepeat2}.
\end{align}
Multiplying \eqref{replace}, \eqref{rinseandrepeat1}, and \eqref{rinseandrepeat2}, and canceling $\prod_{x \in S \setminus B} x_{1}x_{2}x_{3}$ from both sides, we see that \eqref{lemmaoftangents} holds. \qed

\vspace{0.25cm} 

Now we use \lref{segrelemmaoftangents} to show that the matrix $L$ defined in \eqref{entriesofL} does not have full column rank over $\mathbb{F}_{q}$. 

\begin{lemma}
\label{lem:rankL}
Let $S \subset \mathbb{F}_{q}^{k}$ be an arc. If $L$ is the matrix defined in \eqref{entriesofL}, then $L$ does not have full column rank over $\mathbb{F}_{q}$.
\end{lemma}

\noindent \textbf{Proof.} Write $S = \{s_{1}, \ldots, s_{|S|} \}$ in order. We use the ordering of $S$ to write the elements of $A \in \binom{S}{k-2}$, and the elements of $S \setminus A \in \binom{S}{|S|-k+2}$ in order as $A = \{a_{1}, \ldots, a_{k-2} \}$ and $S \setminus A = \{\bar{a}_{1}, \ldots, \bar{a}_{|S|-k+2} \}$. 

Let $L_{(A,A')}$ denote the row of $L$ that is indexed by $(A,A')$. To prove that $L$ does not have full column rank over $\mathbb{F}_{q}$, we will show that the rows in $L$ are spanned by
\begin{equation}
\label{rowspanL}
\mathcal{R} = \{ L_{(\{\bar{a}_{1}\} \cup \{a_{1}, \ldots, a_{k-3}\}, A)} : A \neq \{s_{1}, \ldots, s_{k-2} \} \}. 
\end{equation}
To accomplish this, we must order the rows of $L$. First, list the rows of $\mathcal{R}$ and then list the remaining rows in lex order. We will show that each row of $L$ that is not in $\mathcal{R}$ can be written as a linear combination of two rows of $L$ that precede it. Hence, by induction, every row of $L$ can be written as a linear combination of rows in $\mathcal{R}$.

Let $L_{(A,A')}$ be a row of $L$ that is not in $\mathcal{R}$. We distinguish two cases. 

\vspace{0.25cm}

\noindent \textbf{Case 1:} There exists $s \in S \setminus (A \cap A')$ such that $s$ precedes $A \setminus A'$ and $A' \setminus A$ in the ordering of $S$. 

\vspace{0.25cm}

\noindent Let $\hat{A} = \{s\} \cup (A \cap A')$ and note that $\hat{A} < A < A'$ in lex order. Also observe that $\tau(A, \hat{A} \cup A) = \tau(A', \hat{A} \cup A')$, $\tau(\hat{A}, \hat{A} \cup A') = \tau(A, A \cup A')$, and $\tau(A', A \cup A') = \tau(\hat{A}, \hat{A} \cup A) + 1$.
Let $t$ be defined by \eqref{defoft} and define $w_{1} = (\tau(A, \hat{A} \cup A) + \tau(A, A \cup A') + 1)(t+1) + (t+2)$ and $w_{2} = (\tau(A', \hat{A} \cup A') + \tau(A', A \cup A'))(t+1)$. Applying
\lref{segrelemmaoftangents} with $D = A \cap A'$, $u = A \setminus A'$, $v = A' \setminus A$ and $w=s$ 
implies that $L_{(A,A')}$ is a linear combination of $L_{(\hat{A},A)}$ and $L_{(\hat{A},A')}$:
\begin{equation}
\label{lincomb1}
L_{(A,A')} = (-1)^{w_{1}} \frac{f_{A,S}(A' \setminus A)}{f_{A,S}(s)}L_{(\hat{A},A)} + (-1)^{w_{2}} \frac{f_{A',S}(A \setminus A')}{f_{A',S}(s)}L_{(\hat{A},A')}.
\end{equation}
If $A \cap A' = \{a_{1}, \ldots, a_{k-3}\}$ and $s = \bar{a}_{1}$, then $L_{(\hat{A},A)} \in \mathcal{R}$; otherwise the rows $L_{(\hat{A},A)}$ and $L_{(\hat{A}, A')}$ precede $L_{(A,A')}$.

\vspace{0.25cm}

\noindent \textbf{Case 2: } There does not exist $s \in S \setminus (A \cap A')$ such that $s$ precedes $A \setminus A'$ and $A' \setminus A$ in the ordering of $S$. 

\vspace{0.25cm}

\noindent Write $D = A \cap A' = \{d_{1}, \ldots, d_{k-3}\}$ using the ordering of $S$. Observe that $A \setminus A'$ precedes $d_{k-3}$ in the ordering of $S$; otherwise $D = \{s_{1}, \ldots, s_{k-3}\}$ and $A \setminus A' = s_{k-2}$, which would imply $L_{(A, A')} \in \mathcal{R}$. Let $\hat{A} = D \setminus \{d_{k-3}\} \cup \{A \setminus A'\} \cup \{A' \setminus A\}$ and note that $\hat{A} < A < A'$ in lex order because $A \setminus A'$ precedes $d_{k-3}$ in the ordering of $S$. Since the union of any two of $A$, $A'$, and $\hat{A}$ equals the union of all three, we have $L_{(A,A')} = -L_{(\hat{A},A)} + L_{(\hat{A},A')}$. Also, the rows $L_{(\hat{A},A)}$ and $L_{(\hat{A},A')}$ precede $L_{(A,A')}$. \qed

\vspace{0.25cm}

We now prove \lref{crux} and thus complete the proofs of \tref{mainnew} and \tref{similarr}.

\vspace{0.25cm}

\noindent \textbf{Proof of \lref{crux}} By \lref{rankL}, the matrix $L$ defined in \eqref{entriesofL} does not have full column rank over $\mathbb{F}_{q}$, so there exists a nonzero vector $\vec{\alpha}$ in the nullspace of $L$. The coordinates $\alpha_{A}$ of $\vec{\alpha}$ satisfy \eqref{rowentriesequal} and we now show that they are all nonzero. Suppose, for a contradiction, that there exists $\hat{A} \in \binom{S}{k-2}$ such that $\alpha_{\hat{A}} = 0$. By \eqref{rowentriesequal}, $\alpha_{A'} = 0$ for all $A' \in \binom{S}{k-2}$ such that $\hat{A} \cup A' \in \binom{S}{k-1}$. Repeating this argument, we see that $\alpha_{A} = 0$ for all $A \in \binom{S}{k-2}$, which contradicts that $\vec{\alpha} \neq 0$. Therefore, all coordinates of $\vec{\alpha}$ are nonzero so \lref{crux} holds by \lref{nullL}. \qed  

Let $F$ be the subset consisting of the first $k-2$ elements of $S$. For a subset $A \in \binom{S}{k-2}$, let $D = A \cap F$, let $A \setminus F = \{x_{1}, \ldots, x_{r} \}$, let $F \setminus A = \{z_{1}, \ldots, z_{r} \}$, and let $s$ be the minimum number of transpositions required to order $(F \cap A, F \setminus A)$ as $F$. Let $t$ be defined by \eqref{defoft}. One can show that an explicit solution for a nonzero vector $\vec{\alpha} \in \mathbb{F}_{q}^{\binom{|S|}{k-2}}$ in the nullspace of $L$ is given by 
\begin{equation}
\label{uglyratio}
\alpha_{A} = (-1)^{ (r+s)(t+1) } \prod_{i=1}^{r} \frac{ f_{D \cup \{z_{r}, \ldots, z_{i}, x_{i-1}, \ldots, x_{1}\},S}(x_{i}) }{ f_{D \cup \{z_{r}, \ldots, z_{i+1}, x_{i}, \ldots, x_{1}\},S}(z_{i})},
\end{equation}
which motivates Ball's definition of $\alpha_{A}$ in \cite[Section 3]{Ballextension}.

\section{Proof of \tref{none}}
\label{sec:rankP1}

Let $G \subset \mathbb{F}_{q}^{k}$ be an arc of size $2k-2$. For $C \in \binom{G}{k-1}$, let $e(C) \in \mathbb{F}_{q}^{\binom{2k-2}{k-1}}$ be the $C$-coordinate vector; that is $e(C)_{C'} = 1$ if $C = C'$ and $e(C)_{C'} = 0$ otherwise. To prove that the matrix $M_{G}^{\uparrow 1}$ defined in \eqref{entriesofMGr} has full row rank over $\mathbb{F}_{q}$ when $k \leq 2p-2 \leq q$, we will show that for each $C \in \binom{G}{k-1}$, the $C$-coordinate vector $e(C)$ lies in the column space of $M_{G}^{\uparrow 1}$. 

For a fixed $U \in \binom{G}{1}$, recall that we noted in \sref{ourresults} that the submatrix $M_{G}^{\uparrow 1}(U)$ of $M_{G}^{\uparrow 1}$ equals $D_{U}I_{2k-3}(k-1,k-2)$. Hence, to understand the column space of $M_{G}^{\uparrow 1}$, we must understand the column space of $I_{2k-3}(k-1,k-2)$. By \tref{FWprankformula}, the inclusion matrix $I_{2k-3}(k-1,k-2)$ is invertible over $\mathbb{F}_{q}$ exactly when $k \leq p = \ch(\mathbb{F}_{q})$ so our first goal is to determine a spanning set for the orthogonal space of the column space of $I_{2k-3}(k-1,k-2)$ over $\mathbb{F}_{q}$ when $k > p$. This will allow us to prove that a vector $\vec{y}$ lies in the column space of $I_{2k-3}(k-1,k-2)$ over $\mathbb{F}_{q}$ by showing that $\vec{y}$ is orthogonal to every vector in the spanning set.

\begin{lemma}
\label{lem:colperp}
If $k > p = \ch(\mathbb{F}_{q})$ then, over $\mathbb{F}_{q}$, the nullspace of the inclusion matrix $I_{2k-3}(k+p-2, k-1)$ is the column space of $I_{2k-3}(k-1,k-2)$.
\end{lemma}

\noindent \textbf{Proof.} Over $\mathbb{F}_{q}$, the column space of $I_{2k-3}(k-1,k-2)$ clearly lies in the nullspace of the inclusion matrix $I_{2k-3}(k+p-2,k-1)$ so it suffices to show that the nullity of $I_{2k-3}(k+p-2, k-1)$ equals the rank of $I_{2k-3}(k-1,k-2)$. Observe that the inclusion matrix $I_{2k-3}(k-2, k-p-1)$ equals the inclusion matrix $I_{2k-3}(k+p-2, k-1)^{\top}$ so by \tref{FWprankformula} and Lucas' Theorem \cite{Lucas}, 
\begin{equation}
\label{dimrow}
\nul_{\mathbb{F}_{q}} \, I_{2k-3}(k+p-2, k-1) = \binom{2k-3}{k-2} - \sum_{i \in J} \binom{2k-3}{i} - \binom{2k-3}{i-1},
\end{equation}
where $J = \{ 0 \leq i \leq k-2 : i = k-1 \; \mbox{(mod $p$)} \}$. On the other hand, by \tref{FWprankformula} and Lucas' Theorem,
\begin{equation}
\label{dimcol}
\rank_{\mathbb{F}_{q}} \, I_{2k-3}(k-1,k-2) = \sum_{i \in L} \binom{2k-3}{i} - \binom{2k-3}{i-1},
\end{equation}
where $L = \{ 0 \leq i \leq k-2 : i \neq k-1 \; \mbox{(mod $p$)} \}$. Since \eqref{dimrow} equals \eqref{dimcol}, the lemma follows. \qed

We now define some special vectors in $\mathbb{F}_{q}^{\binom{2k-3}{k-1}}$. We will show in the proof of \tref{none} that variants of these vectors lie in the column space of $M_{G}^{\uparrow 1}$. 

\begin{definition}
\label{def:poweroftwo}
For $0 \leq i \leq k-2$, suppose that $X = \{x_{1}, \ldots, x_{i}\}$, $Y = \{y_{1}, \ldots, y_{i} \}$, and $\Delta = \{y_{i+1}, \ldots, y_{k-1} \}$ are disjoint subsets of $\{1, \ldots, 2k-3\}$. For $\tau \subseteq \{1, \ldots, i\}$, let $X_{\tau} = \{x_{j} : j \in \tau \}$ and let $Y_{\tau} = \{y_{j} : j \in \tau \}$. Define the vector $\vec{v}_{i}(X,Y,\Delta) \in \mathbb{F}_{q}^{\binom{2k-3}{k-1}}$ with coordinates indexed by $\binom{\{1, \ldots, 2k-3\}}{k-1}$ as
\begin{equation}
\label{v}
\vec{v}_{i}(X,Y,\Delta)_{C} = \begin{cases}
                                (-1)^{|\tau|} & \mbox{if $C = X_{\tau} \cup (Y \setminus Y_{\tau}) \cup \Delta$ for $\tau \subseteq \{1, \ldots, i\}$} \\
                                0 & \mbox{otherwise}.
                                \end{cases}
\end{equation}
\end{definition}

We now show that if $k > p = \ch(\mathbb{F}_{q})$, the vector $\vec{v}_{k-p}(X,Y,\Delta)$ defined in \eqref{v} lies in the column space of $I_{2k-3}(k-1,k-2)$ over $\mathbb{F}_{q}$.

\begin{lemma}
\label{lem:base}
If $k > p = \ch(\mathbb{F}_{q})$, then for any choice of $X$, $Y$, and $\Delta$ satisfying the constraints in \dref{poweroftwo}, the vector $\vec{v}_{k-p}(X,Y,\Delta)$ defined in \eqref{v} lies in the column space of $I_{2k-3}(k-1,k-2)$ over $\mathbb{F}_{q}$. 
\end{lemma}

\noindent \textbf{Proof.} By \lref{colperp}, it suffices to show that the vector $\vec{v}_{k-p}(X,Y,\Delta)$ lies in the nullspace of the inclusion matrix $I_{2k-3}(k+p-2, k-1)$ for any choice of $X$, $Y$, and $\Delta$ satisfying the constraints in \dref{poweroftwo}. For $H \in \binom{\{1, \ldots, 2k-3 \}}{k+p-2}$, let $I_{2k-3}(k+p-2, k-1)_{H}$ be the row of the inclusion matrix $I_{2k-3}(k+p-2, k-1)$ corresponding to $H$. We want to show that 
\begin{equation}
\label{dot}
I_{2k-3}(k+p-2, k-1)_{H} \vec{v}_{k-p}(X,Y,\Delta) = 0.
\end{equation}

Define $\overline{H} = \{1, \ldots, 2k-3\} \setminus H$ and define $R(X) = \{ 1 \leq i \leq k-p : x_{i} \in \overline{H} \}$ and $R(Y) = \{ 1 \leq i \leq k-p : y_{i} \in \overline{H} \}$. Moreover, define $\cF$ to be the family of $(k-1)$-subsets $C$ of $\{1, \ldots, 2k-3\}$ such that $I_{2k-3}(k+p-2, k-1)_{(H,C)} \neq 0$ and $\vec{v}_{k-p}(X,Y,\Delta)_{C} \neq 0$. If $H$ does not contain $\Delta$ or if $R(X)$ and $R(Y)$ have nonempty intersection, then $\cF = \emptyset$ and thus \eqref{dot} holds. Otherwise, the elements of $\cF$ are of the form $C = \Delta \cup X_{\tau} \cup (Y \setminus Y_{\tau})$ where $\tau = R(Y) \cup U$ and $U \subseteq \{1, \ldots, k-p \} \setminus (R(X) \cup R(Y))$. Let $W = \{1, \ldots, k-p \} \setminus (R(X) \cup R(Y))$ and observe that $W \neq \emptyset$ because $R(X)$ and $R(Y)$ are disjoint and because $|\overline{H}| = k-p-1$. Since the left hand side of \eqref{dot} equals
\begin{equation}
\label{lastcase}
\sum_{\substack{\tau = R(Y) \cup U \\ U \subseteq W }} (-1)^{|\tau|} = (-1)^{|R(Y)|} \sum_{j=0}^{|W|} \binom{|W|}{j} (-1)^{j} = (-1)^{|R(Y)|}(1-1)^{|W|} = 0,
\end{equation}
the lemma follows. \qed

For a fixed $U \in \binom{G}{1}$, recall that we noted in \sref{ourresults} that the submatrix $M_{G}^{\uparrow 1}(U)$ of $M_{G}^{\uparrow 1}$ equals $D_{U}I_{2k-3}(k-1,k-2)$. 
Since \lref{base} shows that the vector $\vec{v}_{k-p}(X,Y,\Delta)$ lies in the column space of the inclusion matrix $I_{2k-3}(k-1,k-2)$ when $k > p$, we have that for $U \in \binom{G}{1}$, the vector $D_{U} \vec{v}_{k-p}(X,Y,\Delta)$ lies in the column space of the submatrix $M_{G}^{\uparrow 1}(U)$. Padding the vector $D_{U}\vec{v}_{k-p}(X,Y,\Delta)$ with zeroes in the appropriate places thus gives a vector in the column space of $M_{G}^{\uparrow 1}$. To make this precise, we now define variants of the vectors $\vec{v}_{i}(X,Y,\Delta)$ defined in \dref{poweroftwo}.

\begin{definition}
\label{def:poweroftwoP1}
Let $G \subset \mathbb{F}_{q}^{k}$ be an arc of size $2k-2$ and let $U \in \binom{G}{1}$. Let $B$ be the basis of $\mathbb{F}_{q}^{k}$ from \dref{M}. For $0 \leq i \leq k-2$, let $X = \{x_{1}, \ldots, x_{i}\}$, $Y = \{y_{1}, \ldots, y_{i} \}$, and $\Delta = \{y_{i+1}, \ldots, y_{k-1} \}$ be disjoint subsets of $G \setminus U$. For $\tau \subseteq \{1, \ldots, i\}$, let $X_{\tau}$ and $Y_{\tau}$ be defined as in \dref{poweroftwo}. Define the vector $\vec{v}_{i}(U,X,Y,\Delta) \in \mathbb{F}_{q}^{\binom{2k-2}{k-1}}$ with coordinates indexed by $\binom{G}{k-1}$ as
\begin{equation}
\label{vR}
\vec{v}_{i}(U,X,Y,\Delta)_{C} = \begin{cases}
                                (-1)^{|\tau|}\det(U,C)_{B} & \mbox{if $C = X_{\tau} \cup (Y \setminus Y_{\tau}) \cup \Delta$ for $\tau \subseteq \{1, \ldots, i\}$} \\
                                0 & \mbox{otherwise}.
                                \end{cases}
\end{equation}
\end{definition}
Observe that the vector $\vec{v}_{k-p}(U,X,Y,\Delta)$ is the vector $D_{U}\vec{v}_{k-p}(X,Y,\Delta)$ padded with zeroes in all coordinates $C \in \binom{G}{k-1}$ that have nonempty intersection with $U$. Consequently, when $k > p$, the vector $\vec{v}_{k-p}(U,X,Y,\Delta)$ lies in the column space of $M_{G}^{\uparrow 1}$ for any choice of $U \in \binom{G}{1}$ and any choice of $X$, $Y$, and $\Delta$ satisfying the constraints in \dref{poweroftwoP1}.

In the proof of \tref{none}, we show that each of the $C$-coordinate vectors $e(C)$ are linear combinations of the vectors $\vec{v}_{k-p}(U,X,Y,\Delta)$, and hence lie in the column space of $M_{G}^{\uparrow 1}$. To specify the linear combination, we need the following lemma.

\begin{lemma}
\label{lem:linearcomb}
Let $G \subset \mathbb{F}_{q}^{k}$ be an arc of size $2k-2$ and let $B$ be the basis of $\mathbb{F}_{q}^{k}$ fixed in \dref{M}. Let $C \in \binom{G}{k-1}$ and suppose that $\Delta \subset C$. Let $W \subset G \setminus \Delta$ have size $k-|\Delta|$. For any $u \in G$, we have
\begin{equation}
\label{cramerlike}
\sum_{w \in W} \frac{\det(u, W \setminus w, \Delta)_{B}}{\det(w, W \setminus w, \Delta)_{B}} \det(w,C)_{B} = \det(u,C)_{B}. 
\end{equation}
\end{lemma}

\noindent \textbf{Proof.} Since $W \cup \Delta$ is a basis of $\mathbb{F}_{q}^{k}$, we can write $u \in G$ as a unique linear combination of the elements of $W \cup \Delta$. It is easy to see that the coefficient of $w$ in this linear combination is $\det(u, W \setminus w, \Delta)_{B}/\det(w, W \setminus w, \Delta)_{B}$. Since $\Delta \subset C$, \eqref{cramerlike} holds. \qed

\vspace{0.25cm}

The vectors $v_{i}(U,X,Y,\Delta)$ have three nice properties: For fixed $X$, $Y$, and $\Delta$ and $U \in \binom{G \setminus (X \cup Y \cup \Delta)}{1}$, the support of the vector $v_{i}(U,X,Y,\Delta)$ is always the same. Moreover, all the $(k-1)$-subsets $C$ in the support of $v_{i}(U,X,Y,\Delta)$ contain the same fixed set $\Delta$ and have empty intersection with $G \setminus (X \cup Y \cup \Delta)$. Consequently, we can add vectors $v_{i}(U,X,Y,\Delta)$ for different $U \in \binom{G \setminus (X \cup Y \cup \Delta)}{1}$ using \lref{linearcomb} to yield vectors with smaller weight in the column space of $M_{G}^{\uparrow 1}$. Eventually, we conclude that the $C$-coordinate vectors $e(C)$, which have weight one, lie in the column space of $M_{G}^{\uparrow 1}$.

\vspace{0.25cm}

\noindent \textbf{Proof of \tref{none}} The matrix $M_{G}^{\uparrow 1}$ defined in \eqref{entriesofMGr} has full row rank if and only if its column space contains the $C$-coordinate vector $e(C)$ for each $C \in \binom{G}{k-1}$. If $k \leq p$, then the inclusion matrix $I_{2k-3}(k-1,k-2)$ is invertible by \tref{FWprankformula}. For a fixed $U \in \binom{G}{1}$, recall that we noted in \sref{ourresults} that the submatrix $M_{G}^{\uparrow 1}(U)$ of $M_{G}^{\uparrow 1}$ equals $D_{U}I_{2k-3}(k-1,k-2)$. Hence, for each $U \in \binom{G}{1}$, the submatrix $M_{G}^{\uparrow 1}(U)$ is invertible. Thus, the column space of $M_{G}^{\uparrow 1}$ contains the $C$-coordinate vector $e(C)$ for each $C \in \binom{G}{k-1}$.

Now suppose that $p < k \leq 2p-2 \leq q$. Observe that for any $U \in \binom{G}{1}$ and $C \in \binom{G \setminus U}{k-1}$, the $C$-coordinate vector $e(C)$ is a nonzero scalar multiple of the vector $\vec{v}_{0}(U, \emptyset, \emptyset, C)$. We show that the vectors $\vec{v}_{0}(U, \emptyset, \emptyset, C)$ lie in the column space of $M_{G}^{\uparrow 1}$ by proving that for any $0 \leq i \leq k-p$, $U \in \binom{G}{1}$, and $X$, $Y$, and $\Delta$ satisfying the constraints in \dref{poweroftwoP1}, the vector $\vec{v}_{i}(U,X,Y,\Delta)$ defined in \eqref{vR} lies in the column space of $M_{G}^{\uparrow 1}$. 

The proof is by induction on $i$. By \lref{base} and the remarks preceding and following \dref{poweroftwoP1}, the statement is true for the base case $i=k-p$. We assume the statement is true for $i \in \{1, \ldots, k-p\}$ and prove the statement for $i-1 \in \{0, \ldots, k-p-1 \}$. Let $U \in \binom{G}{1}$ and let $X' = \{x_{1}, \ldots, x_{i-1} \}$, $Y' = \{y_{1}, \ldots, y_{i-1} \}$, and $\Delta' = \{y_{i}, \ldots, y_{k-1} \}$ be disjoint subsets of $G \setminus U$. We will show that the vector $\vec{v}_{i-1}(U,X',Y',\Delta')$ lies in the column space of $M_{G}^{\uparrow 1}$.

Define $x_{i} = U$ and let $X = X' \cup \{x_{i}\}$, $Y = Y' \cup \{y_{i}\}$, and $\Delta = \Delta' \setminus \{y_{i}\}$. Write $G$ as the disjoint union $G = X \cup Y \cup \Delta \cup W \cup \Omega$, where $|W|=i+1$ and $|\Omega| = k-2-2i$. Note that $i \leq k-p$ and $k \leq 2p-2$ imply that $|\Omega| \geq 0$. For each $w \in W$, we have that $X$, $Y$, and $\Delta$ are disjoint subsets of $G \setminus w$ satisfying the constraints of \dref{poweroftwoP1}. Consequently, by the induction hypothesis, for each $w \in W$, the vector $\vec{v}_{i}(w,X,Y,\Delta)$ lies in the column space of $M_{G}^{\uparrow 1}$. Moreover, the support of $\vec{v}_{i}(w,X,Y,\Delta)$, denoted $\cS_{w}$, is the same for all $w \in W$,
\begin{equation}
\label{supportvw}
\cS_{w} = \left \{C \in \binom{G}{k-1} : C = \Delta \cup X_{\tau} \cup (Y \setminus Y_{\tau}) \; \mbox{for $\tau \subseteq \{1, \ldots, i\}$} \right \}.
\end{equation}
We now show that the vector $\vec{v}_{i-1}(U,X',Y',\Delta')$ is a linear combination of the vectors $\vec{v}_{i}(w,X,Y,\Delta)$,
\begin{equation}
\label{lincombvr}
\vec{v}_{i-1}(U,X',Y',\Delta') = \sum_{w \in W} \frac{ \det(U, W \setminus w, \Delta)_{B} }{ \det(w, W \setminus w, \Delta)_{B} } \vec{v}_{i}(w,X,Y,\Delta).
\end{equation}
Observe that the support of $\vec{v}_{i-1}(U, X', Y', \Delta')$, denoted $\cS_{U}$, is a subset of the support $\cS_{w}$ from \eqref{supportvw},
\begin{equation}
\label{supportsubset}
\cS_{U} = \left \{C \in \binom{G}{k-1} : C = \Delta \cup X_{\tau} \cup (Y \setminus Y_{\tau}) \; \mbox{for $\tau \subseteq \{1, \ldots, i-1\}$} \right \}.
\end{equation}
Consequently, to prove \eqref{lincombvr}, we must show that the $C$-coordinates of the left and right hand sides of \eqref{lincombvr} are equal,
\begin{equation}
\label{Ccase}
\sum_{w \in W} \frac{ \det(U, W \setminus w, \Delta)_{B} }{ \det(w, W \setminus w, \Delta)_{B} } \vec{v}_{i}(w,X,Y,\Delta)_{C}
= \begin{cases}
   (-1)^{|\tau|} \det(U,C)_{B} & \mbox{if $C \in \cS_{w} \cap \cS_{U}$} \\
   0 & \mbox{if $C \in \cS_{w} \setminus \cS_{U}$}.
  \end{cases}
\end{equation}
We see that \eqref{Ccase} follows from \lref{linearcomb} because if $C \in \cS_{w}$ then $\Delta \subset C$ and $W \subset G \setminus \Delta$ has size $k - |\Delta|$ so the left hand side of \eqref{Ccase} equals
\begin{equation}
\label{almostdone}
\sum_{w \in W} \frac{ \det(U, W \setminus w, \Delta)_{B} }{ \det(w, W \setminus w, \Delta)_{B} } (-1)^{|\tau|} \det(w,C)_{B} = (-1)^{|\tau|} \det(U,C)_{B}.
\end{equation}
If $C \in \cS_{w} \setminus \cS_{U}$, then $i \in \tau$ which implies that $U \in C$ and hence $\det(U,C)_{B}=0$. \qed

\section{Classification}
\label{sec:classification}

To prove \tref{classificationgeneral}, we first state a sufficient condition for an arc $S \subset \mathbb{F}_{q}^{k}$ of size $q+1$ to be linearly equivalent to the normal rational curve $\cR_{k}$.

\begin{lemma}[Roth--Lempel \cite{RothLempel}] 
\label{lem:sufficientforRS}
Suppose that $S \subset \mathbb{F}_{q}^{k}$ is an arc of size $q+1$ and let $B = (e_{1}, \ldots, e_{k}) \subset S$ be a basis of $\mathbb{F}_{q}^{k}$. For $x \in S \setminus B$, let $x = (x_{1}, \ldots, x_{k})$ be the coordinates of $x$ when written with respect to $B$. Let $W_{S,B}$ be a matrix whose columns are the vectors $(x_{1}^{-1}, \ldots, x_{k}^{-1})^{\top}$ for $x \in S \setminus B$. If $\rank W_{S,B} = 2$, then $S$ is linearly equivalent to the normal rational curve $\cR_{k} \subset \mathbb{F}_{q}^{k}$. 
\end{lemma}

Suppose that $S \subset \mathbb{F}_{q}^{k}$ is an arc of size $q+1$ and that there exists a nonnegative integer $n$ for which the hypothesis of \tref{classificationgeneral} is satisfied. Moreover, let $B = (e_{1}, \ldots, e_{k}) \subset S$ be a basis of $\mathbb{F}_{q}^{k}$. To prove that the matrix $W_{S,B}$ defined in \lref{sufficientforRS} has $\rank W_{S,B} = 2$, we will show that any three columns of $W_{S,B}$ are linearly dependent. Given three columns of $W_{S,B}$, we will show they are dependent by constructing a $(k-2) \times k$ matrix $Z$ with $\rank Z = k-2$ so that the three columns of $W_{S,B}$ lie in the nullspace of $Z$. In other words, we want to find $k-2$ independent vectors in $\mathbb{F}_{q}^{k}$ that are orthogonal to each of the three given columns of $W_{S,B}$. Using the notation of \lref{sufficientforRS}, observe that for $x \in S \setminus B$ and $1 \leq j \leq k$, we have $x_{j}^{-1} = (-1)^{j+1} \det(x, B \setminus \{e_{j}\})^{-1}_{B}$. The expression $\det(x, B \setminus \{e_{j}\})_{B}$ has appeared before, for example in \eqref{entriesofF3}, which suggests how to find the required vectors. The following lemma makes this intuition precise.

\begin{lemma}
\label{lem:orthogonal}
Suppose that $0 \leq n \leq q-2k$ and that for every arc $G \subset \mathbb{F}_{q}^{k}$ of size $2k-2+n$, the column space of the matrix $H_{G}^{\uparrow n}$ defined in \dref{ratiodiagonal} contains a vector $v \in \mathbb{F}_{q}^{\binom{2k-2+n}{k-1}}$ such that $v_{i} = 1$ if $i \in \{1, \ldots, k \}$ and $v_{i} = 0$ otherwise. If $S \subset \mathbb{F}_{q}^{k}$ is an arc of size $q+1$ and $B = (e_{1}, \ldots, e_{k}) \subset S$ is a basis of $\mathbb{F}_{q}^{k}$, then there exist nonzero constants $c_{1}, \ldots, c_{k} \in \mathbb{F}_{q}$ such that for any $(k-2)$-subset $A \subset S \setminus B$, we have 
\begin{equation}
\label{weird}
\sum_{j=1}^{k} (-1)^{(j+1)(k-1)}c_{j} \prod_{y \in A} y_{j}^{-1} = 0,
\end{equation}
where $y = (y_{1}, \ldots, y_{k})$ is written with respect to the basis $B$. 
\end{lemma}

\noindent \textbf{Proof.}  Let $A \subset S \setminus B$ be a subset of size $k-2$ and let $\hat{L} \subset S \setminus (B \cup A)$ be a subset of size $|\hat{L}| = n$. Define an arc $G$ and its ordering by $G = (B, A, \hat{L})$. Reorder the arc $S$ so that $G$ is the first $2k-2+n$ vectors of $S$.

Since $|S|=q+1$, we have $t=k-2$, where $t$ is defined by \eqref{defoft}. Observe that $|G|=t+k+n$ and that $|S \setminus G| \geq 1$ since $0 \leq n \leq q-2k$. Since the arc $G \subset S$ satisfies the hypotheses of \cref{extension}, we have that $\vec{1}P_{G}^{\uparrow n} = \vec{0}$. By \tref{similarr}, there exist invertible diagonal matrices $D_{1}$ and $D_{2}$ such that $D_{1}P_{G}^{\uparrow n}D_{2} = M_{G}^{\uparrow n}$. Recalling \dref{ratiodiagonal}, we have
\begin{equation}
\label{looktoinclusionr}
\vec{0} = \vec{1}P_{G}^{\uparrow n} = \vec{1} (J_{G}^{\uparrow n} D_{1})^{-1} (J_{G}^{\uparrow n} M_{G}^{\uparrow n}) D_{2}^{-1} \; \; \; \;  \mbox{so} \; \; \; \; \vec{0} = \vec{1} (J_{G}^{\uparrow n} D_{1})^{-1}H_{G}^{\uparrow n}.
\end{equation}

Recalling that $D_{1} = F_{1}F_{3}$ where $F_{1}$ from \lref{cruxr} is defined by the matrix $E_{1}$ in \lref{nullL} and $F_{3}$ is defined by \eqref{entriesofF3}, we see that the $C$-coordinate of $\vec{1} (J_{G}^{\uparrow n} D_{1})^{-1}$ is
\begin{equation}
\label{suggestivecoordr}
(\vec {1} (J_{G}^{\uparrow n} D_{1})^{-1})_{C} = \alpha_{C} \prod_{y \in G \setminus (C \cup L_{C})} \det(y,C)^{-1}_{B}, 
\end{equation}
where $L_{C}$ is the last $n$-subset of $\binom{G \setminus C}{n}$ in colex order.

Note that $L_{C} = \hat{L}$ for all $C \in \binom{B}{k-1}$ since $B$ is the first $k$ elements of $G$ so
\begin{equation}
\label{evenmoresuggestiver}
(\vec {1} (J_{G}^{\uparrow n} D_{1})^{-1})_{B \setminus \{e_{j}\}} = \alpha_{B \setminus \{e_{j}\}} (-1)^{(j+1)(k-1)} \prod_{y \in A} y_{j}^{-1}
\end{equation}
since $y_{j}^{-1} = (-1)^{j+1} \det(y, B \setminus \{e_{j}\})^{-1}_{B}$.

By assumption, the column space of the matrix $H_{G}^{\uparrow n}$ contains a vector $v \in \mathbb{F}_{q}^{\binom{2k-2+n}{k-1}}$ such that $v_{i} = 1$ if $i \in \{1, \ldots, k \}$ and $v_{i} = 0$ otherwise. Since the rows of the matrix $H_{G}^{\uparrow n}$ are in colex order, by \eqref{looktoinclusionr} and \eqref{evenmoresuggestiver},
\begin{equation}
\label{dependsonA}
0 = \langle \vec{1}(J_{G}^{\uparrow n}D_{1})^{-1}, v \rangle = \sum_{j=1}^{k} (-1)^{(j+1)(k-1)} \alpha_{B \setminus \{e_{j}\}} \prod_{y \in A} y_{j}^{-1}. 
\end{equation}
Hence, \lref{orthogonal} follows by setting $c_{j} = \alpha_{B \setminus \{e_{j}\}}$. \qed

\vspace{0.25cm}

We now prove \tref{classificationgeneral}.

\vspace{0.25cm}

\noindent \textbf{Proof of \tref{classificationgeneral}} Let $S \subset \mathbb{F}_{q}^{k}$ be an arc of size $q+1$ and let $B = (e_{1}, \ldots, e_{k}) \subset S$ be a basis of $\mathbb{F}_{q}^{k}$. Since $S$ is an arc, the matrix $W_{S,B}$ defined in \lref{sufficientforRS} has $\rank W_{S,B} \geq 2$ because if a column of $W_{S,B}$ is a multiple of another column of $W_{S,B}$ then two vectors in $S$ are linearly dependent. To prove $\rank W_{S,B} \leq 2$, we will show that any three columns of $W_{S,B}$ are linearly dependent. Let $w,x,z \in S \setminus B$. We will show that there exists a $(k-2) \times k$ matrix $Z$ with $\rank Z = k-2$ such that the columns of $W_{S,B}$ corresponding to $w,x,z \in S \setminus B$ are in the nullspace of $Z$. As $\nul Z = 2$, this proves that the columns of $W_{S,B}$ corresponding to $w,x,z \in S \setminus B$ are linearly dependent.

To construct $Z$, first choose a $(k-2)$-subset $A \subseteq S \setminus (B \cup \{w,x,z\})$, which is possible since $0 \leq n \leq q-2k$. Write $A = \{a_{1}, \ldots, a_{k-2} \}$ and define $A_{i} = A \setminus \{a_{i} \} \cup \{w\}$. By \lref{orthogonal} applied to $A_{i}$ for $1 \leq i \leq k$ we have
\begin{equation}
\label{Ai}
\sum_{j=1}^{k} \left( (-1)^{(j+1)(k-1)} c_{j} \prod_{y \in A \setminus \{a_{i}\}} y_{j}^{-1} \right) w_{j}^{-1} = 0.
\end{equation}
Defining $Z$ to be the $(k-2) \times k$ matrix with $(i,j)$-entry
\begin{equation}
\label{entriesofZ}
Z(i,j) = (-1)^{(j+1)(k-1)}c_{j} \prod_{y \in A \setminus \{a_{i}\}} y_{j}^{-1},
\end{equation}
we see that \eqref{Ai} implies that the column of $W_{S,B}$ corresponding to $w$ lies in the nullspace of $Z$. Repeating the argument above, we similarly have that the columns of $W_{S,B}$ corresponding to $x$ and $z$ lie in the nullspace of $Z$ as well.

To complete the proof, we must show that $\rank Z = k-2$. Multiplying the $j^{\mathrm{th}}$ column of $Z$ by $(-1)^{(j+1)(k-1)}c_{j}^{-1} \prod_{y \in A} y_{j}$ gives a $(k-2) \times k$ matrix $\overline{Z}$ whose rows are $a_{1}, \ldots, a_{k-2}$. Since $a_{1}, \ldots, a_{k-2}$ are linearly independent vectors, $k-2 = \rank \overline{Z} = \rank Z$. \qed

\vspace{0.25cm}

Finally we prove \tref{classificationzero}.

\vspace{0.25cm}

\noindent \textbf{Proof of \tref{classificationzero}} We first show that if $k \leq p = \ch(\mathbb{F}_{q})$ and $G \subset \mathbb{F}_{q}^{k}$ is an arc of size $2k-2$, then the column space of the matrix $H_{G}^{\uparrow 0}$ contains a vector $v \in \mathbb{F}_{q}^{\binom{2k-2}{k-1}}$ such that $v_{i} = 1$ if $i \in \{1, \ldots, k \}$ and $v_{i} = 0$ otherwise. First note that the matrix $H_{G}^{\uparrow 0}$ equals the inclusion matrix $I_{2k-2}(k-1,k-2)$.

Let $B = \{1, \ldots, k \}$. For each subset $A \in \binom{\{1, \ldots, 2k-2\}}{k-2}$, define $l_{A} = |A \cap B|$ and $\beta_{A} = (-1)^{l_{A}} l_{A}! (k-2-l_{A})!$. For each subset $C \in \binom{\{1, \ldots, 2k-2\}}{k-1}$, define $r_{C} = |C \cap B|$. Define $\vec{\beta} \in \mathbb{F}_{q}^{\binom{2k-2}{k-2}}$ to be a vector with coordinates indexed by $\binom{\{1, \ldots, 2k-2\}}{k-2}$ and entries $\vec{\beta}_{A} = \beta_{A}$. Let $\vec{w} = I_{2k-2}(k-1,k-2) \vec{\beta}$. 

Consider the $C$-coordinate of $\vec{w}$. If $C \not \subset B$, then there are $k-1-r_{C}$ subsets $A \in \binom{C}{k-2}$ such that $l_{A} = r_{C}$. The remaining $r_{C}$ subsets $A \in \binom{C}{k-2}$ satisfy $l_{A} = r_{C} - 1$. Consequently,
\begin{equation}
\label{CnotinB}
\vec{w}_{C} = \sum_{A \in \binom{C}{k-2}} \beta_{A} = (k-1-r_{C}) \cdot (-1)^{r_{C}} r_{C}! (k-2-r_{C})! + r_{C} \cdot (-1)^{r_{C}-1} (r_{C}-1)! (k-1-r_{C})! = 0.
\end{equation}
On the other hand, if $C \subset B$, then all $A \in \binom{C}{k-2}$ satisfy $l_{A} = k-2$ so
\begin{equation}
\label{CinB}
\vec{w}_{C} = \sum_{A \in \binom{C}{k-2}} \beta_{A} = (k-1) \cdot (-1)^{k-2}(k-2)!0! = (-1)^{k-2}(k-1)!,
\end{equation}
which is nonzero since $k \leq p$. The first part of \tref{classificationzero} is proved by setting $v = ((-1)^{k-2}/(k-1)!) \vec{w}$ since the rows of $H_{G}^{\uparrow 0}$ are in colex order.

By \tref{classificationgeneral} if $k \leq \min\{p, q/2\}$, the normal rational curve $\cR_{k}$ is the unique arc in $\mathbb{F}_{q}^{k}$ of size $q+1$. By the well-known principle of duality, this implies that if $k \leq p$ and $k \neq (q+1)/2$, then the normal rational curve $\cR_{k}$ is the unique arc in $\mathbb{F}_{q}^{k}$ of size $q+1$. \qed

\medbreak
\noindent{\sc Acknowledgement:} This research was supported in part by the Institute for Pure and Applied Mathematics and by the Institute for Mathematics and its Applications with funds provided by the National Science Foundation and by the Taylor Family Fund. The author is grateful to Universitat Polit\'{e}cnica de Catalunya for hosting her during two research visits. The author thanks Simeon Ball, Abdul Basit, Jan De Buele, Ben Lund, Jeff Kahn, and Nathan Kaplan for many interesting discussions. The author also thanks Simeon Ball and Jan De Buele for computationally verifying \conjref{PErrank} for many arcs.

\end{document}